\documentclass[12pt]{amsart}

\usepackage{cancel}
\usepackage{latexsym, amssymb, amsmath,esint}
\usepackage{soul}
\usepackage{amsfonts, graphicx}
\usepackage[dvipsnames]{xcolor}

\usepackage{hyperref}
\usepackage[alphabetic]{amsrefs}

\usepackage{wasysym,stmaryrd}

\usepackage{mathtools}

\mathtoolsset{showonlyrefs}
\newcommand{\be}{\begin{equation}}
\newcommand{\ee}{\end{equation}}
\newcommand{\beq}{\begin{eqnarray}}
\newcommand{\eeq}{\end{eqnarray}}

\usepackage{amsthm}

\newtheorem{thm}{Theorem}[section]

\newtheorem{lma}[thm]{Lemma}
\newtheorem{prop}[thm]{Proposition}

\theoremstyle{remark}
\newtheorem{rem}[thm]{Remark}
\numberwithin{equation}{section}

\newtheorem{claim}{Claim}[thm]

\def\be{\begin{equation}}
\def\ee{\end{equation}}
\def\bee{\begin{equation*}}
\def\eee{\end{equation*}}

\def\lf{\left}
\def\ri{\right}

\def\Ric{\text{\rm Ric}}
\def\Rm{\text{\rm Rm}}

\def\p{\partial}

\def\rheat{\lf(\frac{\p}{\p t}-\Delta_{g(t)}\ri)}
\def\tr{\operatorname{tr}}

\def\e{\varepsilon}

\def\a{{\alpha}}
\def\b{{\beta}}

\def\R{\mathbb{R}}

\begin{document}

\title[]
{Uniqueness of Ricci flows from singular initial data}

\author{Laura Flower}
\address[Laura Flower]{Mathematics Institute, University of Warwick, Coventry, CV4 7AL, UK.}
\email{laura.flower@warwick.ac.uk}

 \author{Man-Chun Lee}
\address[Man-Chun Lee]{Department of Mathematics, The Chinese University of Hong Kong, Shatin, Hong Kong, China}
\email{mclee@math.cuhk.edu.hk}

\author{Felix Schulze}
\address[Felix Schulze]{Mathematics Institute, University of Warwick, Coventry, CV4 7AL, UK.}
\email{felix.schulze@warwick.ac.uk}

\renewcommand{\subjclassname}{\textup{2020} Mathematics Subject Classification}

\date{\today}

\begin{abstract} We establish uniqueness of Ricci flows with non-negative sectional curvatures and scaling invariant curvature bounds, assuming that the induced distance functions converge to the same initial distance function up to isometry. In the compact case, restricting the size of the initial singular set, we also establish uniqueness under suitable weaker conditions. 
\end{abstract}

\maketitle

\markboth{Laura Flower, Man-Chun Lee, Felix Schulze}{Uniqueness of Ricci flows from singular initial data}

\section{Introduction} \label{sec:introduction}

Let $(M^n, g_0)$ be a smooth, complete Riemannian manifold. A Ricci flow starting at $g_0$ is a smooth one parameter family of metrics $(g(t))_{t\in [0,T)}$ such that
\begin{equation}
    \partial_tg(t) = -2\Ric(g(t)),\quad g(0)=g_0\, .
\end{equation}
Introduced by Hamilton \cite{Hamilton1982}, it constitutes a nonlinear heat equation for the metric and has proven to be a highly successful tool in differential geometry and topology. 

While many of these results rely on understanding the singularity formation and long-time behaviour of the flow, Ricci flow has also emerged as an important tool to construct smoothings of initial singular spaces while preserving suitable lower curvature bounds and geometric control. Important examples include the proof of the Anderson--Cheeger--Colding--Tian conjecture on the structure of volume non-collapsed 3-dimensional Ricci limit spaces by Simon--Topping \cite{SimonTopping2021, SimonTopping2022} (building on the work of Hochard \cite{Hochard2016}) and the proof of Petrunin's smoothing conjecture for non-negatively curved  polyhedral spaces by Lebedeva--Matveev--Petrunin--Shevchishin in dimension three \cite{LebMatPetShev2015} and recently in all dimensions by Bamler--Cabezas-Rivas \cite{BamlerCabezasRivas2026}. Further examples of Ricci flows and Ricci--DeTurck flows constructed starting from singular initial spaces include \cites{
BurkhardtGuim2019, ChuLee,GiesenTopping2011, KochLamm2012, LammSimon2023,Lai2017,LeeTam,LeeLiu2026, Richard2018, Simon2002, Yin2010}. 

To date only little is known about the uniqueness of Ricci flows starting from singular initial data. Our first theorem establishes uniqueness of Ricci flows with non-negative sectional curvature and scaling invariant curvature estimates. 

\begin{thm}\label{thm:compact-1}
Let $g(t)$ and $\tilde g(t)$,  $t\in (0,T]$, be complete smooth solutions to the Ricci flow on the manifolds $M$ and $\tilde M$, respectively. Suppose that there exists $\a>0$ such that 
\begin{enumerate}
\item[(i)] $\mathrm{sec}(g(t)),\mathrm{sec}(\tilde g(t))\geq 0$;
\item[(ii)] $\sup_M|\Rm(g(t))|+\sup_{\tilde M}|\Rm(\tilde g(t))|\leq \a t^{-1}$
\end{enumerate} 
for $t\in (0,T]$. Denote 
\begin{equation}
d_0:=\lim_{t\to 0^+} d_{g(t)},\;\;\tilde d_0:=\lim_{t\to0^+} d_{\tilde g(t)}
\end{equation}
to be the limiting initial distances on $M$ and $\tilde M$ respectively. Assume $\Phi:(M,d_0)\to (\tilde M,\tilde d_0)$ is an isometry. Then $\Phi$ is a smooth diffeomorphism that satisfies $\Phi^*\tilde g(t)=g(t)$ for all   $t\in (0,T]$.
\end{thm}

\begin{rem} (a) Note that we make no further assumption on the structure of the initial data in the above statement.\\
(b) The existence of the limiting initial distances and uniform attainment follows from Lemma \ref{lma:d_0}.\\
(c) Note that for $n\geq4$ non-negative sectional curvature is in general not preserved under Ricci flow. Nevertheless, the conditions of non-negative curvature operator or non-negative complex sectional curvature are preserved and both imply non-negative sectional curvature.\\
(d) By work of Deruelle \cite{Deruelle2016} expanding gradient Ricci solitons with non-negative curvature operator coming out of a metric cone with link $(\mathbb{S}^{n-1},\gamma)$, where $\gamma$ is a smooth metric on $\mathbb{S}^{n-1}$ with curvature operator greater or equal than one, exist and are unique. The above theorem gives an alternative proof of uniqueness without assuming smoothness of the link.
\end{rem}

Our second theorem allows to relax the lower bound on the curvature by restricting the size of the initial singular set. To state it we consider the cone $\mathrm{C}_{\mathrm{PIC2}}$ of weakly $\mathrm{PIC2}$ algebraic curvature operators (also known as non-negative complex sectional curvature) and $\mathrm{Id}$ the identity curvature operator with sectional curvatures equal to one (see \cite{Topping2024} for further exposition). 
\begin{thm}\label{thm:compact-2}
Let $M^n$ be a compact smooth manifold with distance $d_0$ that generates the same topology as $M$. Let $\Sigma\Subset M$ be closed and $g_0$ be a smooth metric on $M_{reg}:=M\setminus \Sigma$ which is locally compatible with  $d_0$, such that 
\begin{enumerate}
\item[(i)] $\mathcal{H}^n_{d_0}(\Sigma)=0$;
\item[(ii)]  $\mathrm{Vol}_{g_0}(M_{reg})<+\infty$; 
\item[(iii)] there exists an open set $U\supset \Sigma$ such that $\mathrm{Rm}_{g_0}\in \mathrm{C}_{\mathrm{PIC2}}$ on $U\setminus \Sigma$.
\end{enumerate}
Suppose $g(t)$ and $\tilde g(t)$ are Ricci flows on $M\times (0,T]$ such that 
\begin{enumerate}
\item[(a)] $d_{g(t)},d_{\tilde g(t)}\to d_0$ uniformly on $M\times M$ as $t\searrow 0$;
\item[(b)] $g(t),\tilde g(t)\to g_0$ in $C^\infty_{loc}(M_{reg})$ as $t\searrow 0$;
\item[(c)] $\sup_M|\Rm(g(t))|+|\Rm(\tilde g(t))|\leq \a t^{-1}$ for some $\a>0$;
\item[(d)] $\Rm(g(t))+\mathrm{Id}$, $ \Rm(\tilde g(t))+\mathrm{Id}\in \mathrm{C}_{\mathrm{PIC2}}$. 
\end{enumerate}
Then $g(t)=\tilde g(t)$ on $M\times (0,T]$.
\end{thm}

\begin{rem} (a) Note that $\Rm(g)+\mathrm{Id}\in \mathrm{C}_{\mathrm{PIC2}}$ implies that $\Ric(g)$ is uniformly bounded from below and thus Lemma \ref{lma:d_0} together with (b) -- (d) gives the existence of a limiting initial distance for both flows.\\
(b) Considering the flows from compact polyhedral spaces constructed in \cite{LebMatPetShev2015} and \cite{BamlerCabezasRivas2026}, Theorem \ref{thm:compact-2} gives uniqueness in the class of flows satisfying conditions (a) -- (d).
\end{rem}

In the compact case, uniqueness of Ricci flows for smooth initial metrics was first established by Hamilton \cite{Hamilton1982} and later simplified by DeTurck \cite{DeTurck1983}. In the non-compact case for solutions with bounded curvature initially and uniformly throughout the evolution, uniqueness was resolved by Chen--Zhu \cite{ChenZhu2006}, see also the proof via energy methods by Kotschwar \cite{Kotschwar2014}. More recently, the second named author established uniqueness in the non-compact case for potentially unbounded inital curvature together with scaling invariant curvature bounds in \cite{Lee2025}. 

For non-smooth initial data in dimension two, Topping and Yin obtained uniqueness for Ricci flows starting from non-atomic Radon measures \cite{ToppingYin2024}. In higher dimensions uniqueness for initial limit spaces with two-sided sectional curvature bounds was proven by Corro, Moreno, and Zarei in \cite{CorroMorenoZarei2025}, provided
the flows also have uniform two-sided sectional curvature bounds. For uniformly volume non-collapsed flows with scaling invariant curvature bounds und a uniform lower bound on the sectional curvature coming out of Reifenberg spaces, uniqueness was recently established by the first-named author in \cite{Flower}.

\medskip

\noindent {\bf Outline of proof.} We employ a similar strategy as in \cite{ChenZhu2006,Lee2025, Flower}, by constructing a harmonic map heat flow $F:M\times(0,T]\to\tilde{M}$ such that
\begin{equation}
    \partial_tF=\Delta_{g(t),\tilde{g}(t)}F
\end{equation}
with $F(x,t)\to\Phi(x)$ as $t\to0$ (where we take $\Phi=\mathrm{id}$ in the case of Theorem \ref{thm:compact-2}). In particular, $F$ is constructed by starting at a positive time $t>0$, smoothing out the map $\Phi:(M,g(t))\to(\tilde{M},\tilde{g}(t))$ to give derivative bounds, flowing this map via harmonic map heat flow, and using the derivative bounds and parabolic regularity to take a smooth limit as $t\to0$. However, unlike as in \cite{Flower}, we use the standard heat flow to obtain a smoothing of $\Phi$ with suitable estimates. The main further important ingredient is an extension of the estimate by Bamler--Brendle \cite{BamlerBrendle} to control the optimal local Lipschitz constant of the map $F$, while allowing for unbounded curvature. This allows to establish that the limiting map $F(\cdot, t)$ is 1-Lipschitz and homotopy equivalent to the initial isometry for all $t>0$. In the compact case, using that the degree of $F(\cdot, t)$ has to be one, and comparing the volumes yields that $F(\cdot, t)$ is an isometry for all $t>0$ and thus constant in time. In the noncompact case, we show that by evolving spherical Busemann functions via the heat equation along the flow, we are able to construct a suitable weighted volume to extend the above argument to this case as well. Note that this differs from the further strategy in both \cite{Lee2025, Flower} which depended on stability estimates for Ricci-DeTurck flow.

The paper is structured as follows. In \S2, we recall properties of the distance function under Ricci flow, and in \S3, we give some useful evolution equations for harmonic map heat flow. Then in \S4, we perform the smoothing of $\Phi$ using heat flow methods. In \S5, we prove Theorem \ref{thm:compact-1} in the compact case, by constructing the harmonic map heat flow $F$ and comparing the volumes of $(M,g(t))$ and $(\tilde{M},\tilde{g}(t))$. We prove Theorem \ref{thm:compact-2} in \S6, by using (almost) nonnegative curvature together with heat kernel estimates to control $F$ near the singular set $\Sigma$, and smoothness back to time $t=0$ to control $F$ away from $\Sigma$. Again, comparing the volumes of $(M,g(t))$ and $(\tilde{M},\tilde{g}(t))$ gives additional control on $F$, which allows us to use properties of the Ricci-DeTurck flow $F_*g(t)$ to prove Theorem \ref{thm:compact-2}. Finally, in \S7, we prove the noncompact case of Theorem \ref{thm:compact-1}, by evolving the spherical Busemann functions via the heat equation along the flows to construct a suitable weighted volume to extend the previous argument to this case.

\medskip

\noindent{\bf Acknowledgments.} L.~Flower is supported by the Warwick Mathematics Institute Centre for Doctoral Training, and gratefully acknowledges funding from the University of Warwick and the UK Engineering and Physical Sciences Research Council (Grant number: EP/W524645/1).  M.-C.~Lee is supported by Hong Kong RGC grants No. 14300623 and No. 14304225, and an Asian Young Scientist Fellowship. F.~Schulze has received funding by the European Union (European Research Council (ERC), GENREG -- 101200301). Views and opinions expressed are however those of the authors only and do not necessarily reflect those of the European Union or the European Research Council. Neither the European Union nor the granting authority can be held responsible for them.

\medskip

\noindent {\bf AI disclosure.} An initial version of the full proof with a more ad-hoc smoothing in the compact case with isolated conical singularities is due to the authors. The idea of smoothing the initial isometry via the standard heat flow was suggested by an LLM, together with a strategy in case the initial data is given by a metric cone. The idea to evolve the Busemann function by the heat flow coupled to Ricci flow to obtain a suitable weighted volume in the non-compact case is again due to the authors.

The authors retain full responsibility for the correctness of the final work. This article does not contain any AI written text.

\section{Canonical distance induced from Ricci flow}

In this section, we discuss the distance which is canonically induced by the Ricci flow. We consider a Ricci flow solution $g(t),t\in (0,T]$ defined on a smooth manifold $M$. In our setting, the distance $d_{g(t)}$ as $t\to0$ is non-degenerating and defines a distance on $M$. The following is proved by Simon-Topping \cite{SimonTopping2021}, building on the distance estimates of Hamilton \cite{Hamilton1995}.

\begin{lma}\label{lma:d_0}
Suppose $g(t),t\in (0,T]$ is a complete solution to the Ricci flow on a smooth manifold $M$ such that for some $\a,c_0>0$,
\begin{equation}
-\a\leq \Ric(g(t))\leq (n-1)c_0 t^{-1}
\end{equation}
on $M\times (0,T]$. Then the distance metric converges uniformly to a distance metric $d_0$ on $M$ as $t\to 0^+$, so that 
\begin{enumerate}
\item[(a)] for all $x,y\in M$ and $t\in (0,T]$,
\begin{equation}
d_0(x,y)-\b(n)\sqrt{c_0t}\leq d_{g(t)}(x,y)\leq e^{\a t}d_0(x,y),
\end{equation}
\item[(b)] for a fixed $x_0\in M$ and any $r>0$, there is $\gamma(n,c_0,r,T)>0$ such that for all $x,y\in B_{g(T)}(x_0,r)$ and $t\in (0,T]$, 
\begin{equation}
d_{g(t)}(x,y)\geq \gamma \cdot [d_0(x,y)]^{1+2(n-1)c_0}.
\end{equation}
\end{enumerate}
In particular, the distance $d_0$ generates the same topology as $M$.
\end{lma}
\begin{proof}

This follows from \cite[Lemma 3.1]{SimonTopping2021} with arbitrarily large $r>0$, and  uniqueness of the limit of $d_{g(t)}$.
\end{proof}

\section{Evolution equations of harmonic map heat flow}
In this section, we collect some evolution equations related to the harmonic map heat flow. We first start with the evolution equation for the pull-back metric under the harmonic map heat flow.

\begin{lma}\label{lma:evo-harm-eign}
Let $(M,g(t))$ and $(N,\tilde g(t))$ be two solutions to the Ricci flow on $[0,T]$. Suppose $F(t):(M,g(t))\to (N,\tilde g(t))$ is a family of smooth
maps satisfying the harmonic map heat flow: 
\begin{equation}\label{eqn:HMHF}
\partial_t F=\Delta_{g(t),\tilde g(t)}F
\end{equation}
on $M\times [0,T]$. Then the tensor $H:=F^* \tilde g(t)$ satisfies
\begin{equation}
\begin{split}
(\partial_t-\Delta_{g(t)}) H_{ij}&=-2g^{k\ell} F^\a_{ik}F^\b_{j\ell} \tilde g_{\a\b}-R_i^k H_{kj}-R_j^k H_{ik}\\
&\quad -2F^\a_i F^\b_j \tilde R_{\a\b}+2 g^{k\ell} F_i^\a F^\b_j F^\gamma_k F^\delta_\ell \tilde R_{\a \gamma\delta \b}.
\end{split}
\end{equation}
Here the connection on the bundle $(T^*M)^{k+1}\otimes F^{-1}(T\tilde M)$ is induced by the Riemannian connections of $g(t)$ and $\tilde g(t)$.
\end{lma}
\begin{proof}
The statement follows directly from the Ricci identity and the Ricci flow equation.
\end{proof}

For the purpose of constructing a harmonic map heat flow, we later will need to control and localize the optimal Lipschitz constant of the map $F$ allowing for unbounded curvature. Extending the argument of Bamler--Brendle in \cite{BamlerBrendle}, we establish the following sharper evolution inequality for the maximal distortion of $DF$.

\begin{lma}\label{lma:evo-eign}
Along the harmonic map heat flow \eqref{eqn:HMHF}, the function $\varphi$ given by 
\begin{equation}
\varphi(x,t):=\sup\{ F^*\tilde g(v,v): g(v,v)=1\}
\end{equation}
satisfies 
\begin{equation}
\rheat \varphi\leq  \varphi (\varphi-1)_+\left(F^*\mathrm{scal}_{\tilde g}+n(n-1)\sigma \right) +2\varphi(n-1)\sigma
\end{equation}
in the barrier sense, where $\sigma$ denotes the negative part of the sectional curvature lower bound of $\mathrm{Rm}(\tilde g)$ at $F(x,t)$.
\end{lma}
\begin{proof}

At a  point $x_0\in M$ and $ F(x_0,t_0)\in \tilde M$, we pick an orthonormal frame $\{e_i\}_{i=1}^n$ with respect to $g(t_0)$ at $x_0$ and an orthonormal frame  $\{u_\a\}_{\a=1}^m$ with respect to $\tilde g(t_0)$ at $F(x_0,t_0)$  such that $F_i^\a=\lambda_i \delta_{i\a}$ with $\lambda_1^2\geq ...\geq \lambda_n^2\geq 0$ at $(x_0,t_0)$. This can be done by first choosing $\{e_i\}_{i=1}^n$ diagonalizing $F^*\tilde g$ and then choosing $u_i:= dF(e_i)/|dF(e_i)|$ for $1\leq i\leq \mathrm{rank}(dF)$ at $(x_0,t_0)$. In particular, $\varphi(x_0,t_0)=(F^*g)(e_1,e_1)=\lambda_1^2$. We extend $\{e_i\}_{i=1}^n$ smoothly as follows: by parallel transport along radial directions using $g(t_0)$ we extend $\{e_i\}_{i=1}^N$ to an orthonormal
frame $\{E_i(x,t_0)\}_{i=1}^n$ around $(x_0,t_0)$. 
We then evolve the frame in time locally by 
\begin{equation}\label{eqn:bundle-trick}
\partial_t E_i^l =R^l_m E_i^m.
\end{equation}
This yields an orthonormal frame  $\{E_i\}_{i=1}^n$ in spacetime near $(x_0,t_0)$ which satisfies 
\begin{equation}\label{eqn:parallel}
\nabla E_i=\Delta E_i=0,\;\;\text{at}\;(x_0,t_0) 
\end{equation}
for $i=1, \ldots n$.

The function $\ell(x,t):=(F^*g)(E_1,E_1)$ is a locally defined smooth function such that $\ell(x_0,t_0)=\varphi(x_0,t_0)$ and $\ell\leq \varphi$. Furthermore by combining \eqref{eqn:bundle-trick},  \eqref{eqn:parallel} with Lemma~\ref{lma:evo-harm-eign}, we have 
\begin{equation}\label{eqn:evo-ell-1} 
\begin{split}
\rheat\Big|_{(x_0,t_0)}\ell(x,t) 
&=\left[\rheat H\right](E_1,E_1)\\
&\quad\;+2H\left(\Ric(E_1),E_1\right)\\
&=
-2g^{kl} F^\a_{1k}F^\b_{1l} \tilde g_{\a\b} -2F^\a_1 F^\b_1 \tilde R_{\a\b}\\
&\quad\; +2 g^{kl} F_1^\a F^\b_1 F^\gamma_k F^\delta_l \tilde R_{\a \gamma\delta \b}\\
&\leq -2\lambda_1^2\tilde R_{11}+2\lambda_1^2  \sum_{k=2}^n\lambda_k^2 \tilde R_{1kk1}\\
&=2\varphi\cdot  \sum_{k=2}^n(\lambda_k^2-1)\tilde R_{1kk1}.
\end{split}
\end{equation}
Here $\tilde R_{1kk1}:=\tilde R(u_1,u_k,u_k,u_1)$ is the sectional curvature of the plane spanned by $\{u_1,u_k\}$ at $F(x_0,t_0)$.

Since by definition of $\sigma$, $\tilde R_{1kk1}\geq -\sigma$, we observe that 
\begin{equation}\label{eqn:evo-ell-2}
\begin{split}
\sum_{k=2}^n(\lambda_k^2-1)\tilde R_{1kk1}&=\sum_{k=2}^n(\lambda_k^2-1)(\tilde R_{1kk1}+\sigma)-\sigma \sum_{k=2}^n (\lambda_k^2-1)\\
&\leq (\varphi-1)\sum_{k=2}^n (\tilde R_{1kk1}+\sigma)+\sigma \sum_{k=2}^n (1-\lambda_k^2)\\
&\leq (\varphi-1)_+\left(\tilde R_{11}+(n-1)\sigma \right) +(n-1)\sigma\\
&\leq \frac12 (\varphi-1)_+\left(\mathrm{scal}_{\tilde g}+n(n-1)\sigma \right) +(n-1)\sigma
\end{split}
\end{equation}
where we have used \cite[Lemma 3.2]{ChanLee2025} applied to $\widetilde \Rm(g(x_0,t_0))+\sigma(x_0,t_0)\mathrm{Id}$ in the last inequality.  The result follows by combining \eqref{eqn:evo-ell-1} and \eqref{eqn:evo-ell-2}.
\end{proof}

\section{Smoothing an almost distance decreasing map}

In this section, we will discuss how the Ricci flows in Theorem~\ref{thm:compact-1} and \ref{thm:compact-2} can be weakly gauged between each other (almost) initially.

\subsection{Smoothing additive Lipschitz maps to Euclidean space}

One of the major difficulties in smoothing the map $\phi$ in Theorem~\ref{thm:compact-1} and \ref{thm:compact-2} lies in the presence of the additive Lipschitz error appearing in (a) of Lemma~\ref{lma:d_0}. We discuss how this can be overcome using the standard heat flow and Lipschitz approximation. We first consider a map from a manifold to Euclidean space. The argument in this sub-section is static.

We need a version of a gradient estimate.
\begin{lma}\label{lma:parab-boot}
Let $(N^n,g)$ be a complete manifold with $\Ric\geq -k, k\geq 0$. Suppose $F:M\times [0,\mathcal{T}]\to\mathbb{R}$ is a solution to the heat equation $(\partial_\tau-\Delta_{g})F=0$. Then
\begin{equation}
\tau|\nabla F|^2(x,\tau)\leq 100(n+k) \cdot \sup_{B_g(x_0,2)\times [0, \mathcal{T}]} F^2
\end{equation}
for all $(x,\tau)\in B_g(x_0,1)\times [0,1 \wedge \mathcal{T}]$.
\end{lma}
\begin{proof}
By direct computation and the Bochner formula, we have 
\begin{equation}\label{eqn:evo-heat}
\left\{
\begin{array}{ll}
\displaystyle\left(\frac{\partial}{\partial \tau}-\Delta_g \right) F^2=-2|\nabla F|^2;\\[3mm]
\displaystyle\left(\frac{\partial}{\partial \tau}-\Delta_g \right) |\nabla F|^2 \leq -2|\nabla^2 F|^2 +2k|\nabla F|^2
\end{array}
\right.
\end{equation}
for $\tau\in (0,\mathcal{T}]$.

Let $\phi$ be a smooth non-increasing function on $[0,+\infty)$ such that $\phi\equiv 1$ on $[0,1]$, vanishes outside $[0,2]$ and satisfies $|\phi'|^2\leq 4\phi$, $\phi''\geq -10\phi$. We let $\Phi(x,\tau):=\phi(d_g(x,x_0))$ be a cut-off function on $N$ which therefore satisfies
\begin{equation}
\begin{split}
    \Delta_g\Phi &\geq
    -10\Phi-2\sqrt{\Phi}\left((n-1)+\sqrt{(n-1)k}\right) \\
    &\geq -10\Phi -(3n+k)\sqrt{\Phi}
\end{split}
\end{equation}
and
\begin{equation}
    |\nabla\Phi|^2 \leq 4\Phi,
\end{equation}
in the barrier sense, where we have used that Laplacian comparison implies that $\Delta_g d_g(x,x_0)\leq (n-1) + \sqrt{(n-1)k}$ if $d_g(x,x_0)\geq 1$.
We consider the function
\begin{equation}
H:=\tau|\nabla F|^2 \Phi+ \Lambda F^2
\end{equation}
for a large $\Lambda>0$.

Using Calabi's trick, we may assume $\Phi$ to be smooth when applying the maximum principle. On the support of $\Phi$, we have
\begin{equation}
\begin{split}
    \left(\frac{\partial}{\partial \tau}-\Delta_g \right)(|\nabla F|^2 \Phi)
    &= \Phi\left(\frac{\partial}{\partial \tau}-\Delta_g \right)|\nabla F|^2 - |\nabla F|^2 \Delta_g\Phi \\
    &\quad\quad -2\langle\nabla |\nabla F|^2,\nabla \Phi \rangle \\
    &\leq |\nabla F|^2 (2k\Phi +10\Phi +(3n+k)\sqrt{\Phi}) \\
    &\quad\quad -2|\nabla^2F|^2\Phi +4|\nabla^2 F||\nabla F||\nabla \Phi|\\
    &\leq |\nabla F|^2 (2k\Phi +10\Phi +(3n+k)\sqrt{\Phi} +8),
\end{split}
\end{equation}
so if $\tau\in (0,1]$ then $H$ satisfies
\begin{equation}\label{eqn:evo-F}
\begin{split}
\left(\frac{\partial}{\partial \tau}-\Delta_g \right) H
&\leq  -2\Lambda|\nabla F|^2+ |\nabla F|^2 \Phi +\tau \left(\frac{\partial}{\partial \tau}-\Delta_g \right)(|\nabla F|^2 \Phi) \\
&\leq |\nabla F|^2\left(-2\Lambda +100(n+l) \right).
\end{split}
\end{equation}

We take $\Lambda=100(n+l)$ and apply the maximum principle to the function $H$ on $B_g(x_0,2)$. Its maximum can only be attained when $t=0$ or on $\partial B_g(x_0,2)$ thanks to \eqref{eqn:evo-F}. This proves the Lemma.
\end{proof}

We also need the following (sharp) preservation of a Lipschitz bound.
\begin{lma}\label{lma:pres-Lips}
Let $(N,g)$ be a complete manifold with $\Ric\geq -k,k\geq 0$. If $f:N\to\mathbb{R}$ is a $L$-Lipschitz function, then it admits a heat flow $F$ from $f$ on $N\times (0,1)$ such that $F$ is $Le^{k\tau}$-Lipschitz and satisfies
\begin{equation}
    |F(x,\tau)-f(x)|\leq C(n,k)\sqrt{\tau}.
\end{equation}
\end{lma}
\begin{proof} 
By the work of Greene-Wu \cite[Proposition 2.1]{GreeneWu}, there exists a sequence of smooth functions $f_i:N\to\mathbb{R}$ such that $f_i\to f$ uniformly, and $f_i$ is $L$-Lipschitz for all $i$. 

By the estimate of Li-Yau \cite[Corollary 3.1]{LiYau}, the heat kernel $G$ of $(N,g)$ satisfies
\begin{equation}
    G(x,y,\tau) \leq \frac{C_0}{\mathrm{Vol}_g(x,\sqrt{t})^{1/2} \mathrm{Vol}_g(y,\sqrt{t})^{1/2}} \exp\left( -\frac{d_g(x,y)^2}{C_0t} \right).
\end{equation}
By taking
\begin{equation}
F_i(x,\tau):=\int_N G(x,y,\tau) f_i(y)\,d\mathrm{vol}_g,
\end{equation}
each $f_i$ admits a heat flow $F_i$. Further, if we fix $x_0\in N$, then $|f_i(x)|\leq Ld_g(x_0,x) + |f_i(x_0)|$ for every $x\in N$, so by the co-area formula, Stokes' Theorem and volume comparison,
\begin{equation}\label{eqn:bound|F|}
    |F_i(x,\tau)|\leq C_1(d_g(x,y) + |f_i(x_0)|)
\end{equation}
for some $C_1=C_1(n,k,L)$.

Next, we use this to obtain a bound on $|\nabla F_i|$. We use an argument similar to the proof of Lemma~\ref{lma:parab-boot}, instead working on $B_g(x_0,2r)$ for arbitrary $x_0\in N$ and $r>1$, with a suitable choice of cutoff function $\Phi_r$. In particular, with $\phi$ as in the proof of Lemma~\ref{lma:parab-boot}, we take $\Phi_r:=\phi(r^{-1}d_g(x,x_0))$, so that
\begin{equation}
\begin{split}
    \Delta_g\Phi_r &\geq -10r^{-2}\Phi_r -(3n+k)r^{-1}\sqrt{\Phi_r},\\
    |\nabla\Phi_r|^2 &\leq 4r^{-2}\Phi_r,
\end{split}
\end{equation}
so on the support of $\Phi$,
\begin{equation}
\begin{split}
    \left(\frac{\partial}{\partial \tau}-\Delta_g \right) \left(|\nabla F_i|^2\Phi_r\right) &\leq |\nabla F_i|^2(2k\Phi_r + 10r^{-2}\Phi_r + (3n+k)r^{-1}\sqrt{\Phi_r} + 8r^{-2}) \\
    &\leq |\nabla F_i|^2(100(n+k)r^{-2}+100(n+k)\Phi_r).
    \end{split}
\end{equation}
Then the function
\begin{equation}
    H:= e^{-100(n+k)\tau}|\nabla F_i|^2\Phi_r + 50(n+k)r^{-2}F_i^2
\end{equation}
satisfies $\left(\frac{\partial}{\partial \tau}-\Delta_g \right)H\leq0$ on the support of $\Phi$. So by applying the maximum principle to $H$ and using \eqref{eqn:bound|F|}, we obtain
\begin{equation}
    |\nabla F_i|^2 \leq C_2(1+|f_i(x_0)|^2r^{-2})
\end{equation}
on $B_g(x_0,r)$ for some $C_2=C_2(n,k,L)$. In particular, we obtain a uniform bound for $|\nabla F_i|^2$ on $N\times(0,1)$.

Then recalling that $\left(\frac{\partial}{\partial \tau}-\Delta_g \right) |\nabla F_i|^2 \leq 2k|\nabla F_i|^2$ and $|\nabla f_i|^2\leq L^2$, the maximum principle for bounded functions gives
\begin{equation}
|\nabla F_i|^2 \leq L^2e^{2k\tau}
\end{equation}
for all $(x,\tau)\in N\times (0,1)$.

Parabolic regularity theory enables us to pass $i\to+\infty$ and obtain a limiting heat flow $F$ from $f$ on $M\times (0,1)$. The assertion of the time zero asymptotic follows from the parabolic estimate $|\partial_\tau F|\leq C(n,k)\tau^{-1/2}$ and integration.
\end{proof}

We are now ready to prove a smoothing result for additive Lipschitz maps to Euclidean spaces. 
\begin{prop}\label{prop:smoothing-map-Euc}
Let $(N,g)$ be a complete manifold with $\Ric\geq -k,k\geq 0$. Suppose $f:N\to\mathbb{R}^N$ is a bounded Borel map such that 
\begin{equation}
|f(x)-f(y)|\leq L d_g(x,y)+\eta.
\end{equation} 
Then $f$ admits a harmonic map heat flow 
\begin{equation}
\left\{
\begin{array}{ll}
\partial_\tau F=\Delta_{g,g_{euc}}F;\\[2mm]
F(0)=f
\end{array}
\right.
\end{equation}
for $\tau\in (0,\hat{\mathcal{T}}]$ so that for all $(x,\tau)\in M\times (0,\hat{\mathcal{T}} ]$,
\begin{enumerate}
\item[(a)] $\|dF\|_{\infty,op} \leq Le^{k\tau}+10\sqrt{n+k}\, \tau^{-1/2}\eta$;
\item[(b)] $|F(x,\tau)-f(x)|\leq C(n,k)L\sqrt{\tau}+\eta$.
\end{enumerate}
\end{prop}
\begin{proof}
We solve the harmonic map heat flow by component-wise convolution with the heat kernel, as in the proof of Lemma~\ref{lma:pres-Lips}.

Fix any $z\in \mathbb{R}^N$ with $|z|=1$. We consider $f_z:=\langle f,z \rangle$, which is a real-valued function on $N$. Define 
\begin{equation}
h_z(x):=\inf_{y\in N} \{f_z(y)+L d_g(x,y)\},
\end{equation}
which is a $L$-Lipschitz function on $N$ with $f_z-\eta\leq h_z \leq f_z$. We denote $f_z:=h_z+g_z$ so that $0\leq g_z\leq \eta$. Then $F_z=\langle z,F\rangle$ can be re-written as the sum of $H_z$ and $G_z$:
\begin{equation}
\left\{
\begin{array}{ll}
\displaystyle H_z(x,\tau):=\int_N G(x,y,\tau)\, h_z(y)\,d\mathrm{vol}_g;\\[2mm]
\displaystyle G_z(x,s\tau):=\int_N G(x,y,\tau)\, g_z(y)\,d\mathrm{vol}_g,
\end{array}
\right.
\end{equation}
where $G(x,y,\tau)$ is the heat kernel of $(N,g)$. By the estimate of Li-Yau \cite[Corollary 3.1]{LiYau} and the maximum principle, we have $0\leq G_z\leq \eta$ on $M\times (0,+\infty)$ and hence  
Lemma~\ref{lma:parab-boot} implies 
\begin{equation}
|\nabla G_z|\leq  10 \sqrt{n+k}\,\tau^{-1/2} \eta
\end{equation}
on $N\times (0,\hat{\mathcal{T}}]$, while Lemma~\ref{lma:pres-Lips} implies $|\nabla H_z|\leq Le^{k\tau}$ for $\tau\in (0,\hat{\mathcal{T}}]$. Since $z$ is arbitrary, this proves the conclusion (a).

To prove (b), by using the kernel representation component-wise and $0\leq G_z\leq \eta$, we have 
\begin{equation}
\begin{split}
|F(x,\tau)-f(x)|&\leq L\int_N G(x,y,\tau)\,d_g(x,y)\,d\mathrm{vol}_g+\eta\\
&\leq C(n,k)L\sqrt{\tau}+\eta
\end{split}
\end{equation}
where we have used Li-Yau's heat kernel estimate \cite{LiYau}.
\end{proof}

\subsection{Smoothing the initial isometry: the compact case}

We first discuss how to approximate the distance isometry when $M$ is compact. We state and prove the approximation which allows a possible regular set.

\begin{prop}\label{prop:gauge-fixing-cpt-1}
Let $g(t)$ and $\tilde g(t)$, $t\in (0,T]$ be smooth solutions to the Ricci flow on compact manifolds $M^n$ and $\tilde M^n$ respectively such that 
\begin{enumerate}
\item[(i)] $\sup_M|\Rm(g(t))|+\sup_{\tilde M}|\Rm(\tilde g(t))|\leq \a t^{-1}$;
\item[(ii)] $\Ric(g(t)),\Ric(\tilde g(t))\geq -(n-1)$
\end{enumerate}
for $t\in (0,T]$. Denote 
\begin{equation}
d_0:=\lim_{t\to 0^+} d_{g(t)},\;\;\tilde d_0:=\lim_{t\to0^+} d_{\tilde g(t)}
\end{equation}
for the distance functions on $M$ and $\tilde M$ respectively. Suppose the following hold:
\begin{enumerate}
\item[(a)] there exists a distance isometry $\Phi:(M,d_0)\to (\tilde M,\tilde d_0)$;
\item[(b)] there exists a closed set $\Sigma\Subset M$ (which could be the whole of $M$) and smooth metrics $g_0,\tilde g_0$ on $M_{reg}:=M\setminus \Sigma$ and $\tilde M_{reg}:=\Phi(M_{reg})$ respectively;
\item[(c)] $g(t)\to g_0$ in $C^\infty_{loc}(M_{reg})$ as $t\to0$;
\item[(d)] $\tilde g(t)\to \tilde g_0$ in $C^\infty_{loc}(\tilde M_{reg})$ as $t\to0$.
\end{enumerate}
Then there exist $t_i \to 0, s_i\in (0,t_i)$, $\Lambda_i \to 1$ and a sequence of smooth maps $f_i: M\to \tilde M$ such that 
\begin{enumerate}
\item[(I)] $f_i^*\tilde g(t_i)\leq \Lambda_i g(s_i)$ on $M$;
\item[(II)] $\sup_M \tilde d_0(f_i(x),\Phi(x))\to 0$ as $i\to+\infty$;
\item [(III)] For any compact set $\Omega\Subset M_{reg}$, there is $N$ such that $f_i=\Phi$ for $i>N$.
\end{enumerate}
\end{prop}
\begin{proof}
By \cite{CalabiHartman}, $\Phi$ is smooth on $M_{reg}$. We choose a compact exhaustion $\Omega_i$ of $M_{reg}$ such that $\Omega_i\Subset \mathrm{int}(\Omega_{i+1})$ for all $i$. Choose a smooth cutoff function $\chi_i$ on $M$ such that $\chi_i=1$ on a neighbourhood of $\Omega_i$ and $\chi_i$ is compactly supported on $\mathrm{int}(\Omega_{i+1})$. 
 
We fix any $t_i\to 0^+$ and $\e_i\to 0^+$ such that $e^{(n-1)t_i}\leq 1+\e_i/2$. By Nash's embedding Theorem \cite[Theorem 3]{Nash}, there is a smooth embedding 
\begin{equation}
\iota_i: (\tilde M^n, \tilde g(t_i)) \to \mathbb{R}^{N}
\end{equation}
for some $N=N(n)$. We let $\pi_i$ be the nearest-point projection on an open neighbourhood $U_i$ of $\iota_i(\tilde M)\subseteq \mathbb{R}^{N}$ to $\tilde M$ with $\|d\pi_i\|_{op}\leq 1+\e_i$.

We  consider the map 
\begin{equation}
\hat f_i:=\iota_i \circ \Phi : M^n\to\mathbb{R}^{N}.
\end{equation}
By Lemma~\ref{lma:d_0}, it is continuous and satisfies 
\begin{equation}
\begin{split}
|\hat f_i(x)-\hat f_i(y)|&\leq d_{\tilde g(t_i)}\left(\Phi(x),\Phi(y) \right)\\
&\leq  e^{(n-1)t_i}\tilde d_{0}\left(\Phi(x),\Phi(y) \right)\\
&=e^{(n-1)t_i} d_0(x,y)\\
&\leq e^{(n-1)t_i} d_{g(s)}(x,y)+C(n,\a) s^{1/2}
\end{split}
\end{equation}
for all $s\in (0,T]$.

We choose $s_i$ small (relative to $t_i$) and apply Proposition~\ref{prop:smoothing-map-Euc} to $(N,g)=(M,g(s_i))$ to obtain a smooth map $F_i:=F(\tau_i)$ for small $\tau_i:=s_i^{1/2}\in (0,1)$, so that $F_i:(M,g(s_i))\to \mathbb{R}^{N}$ satisfies  
\begin{enumerate}
\item[(A)] $\|dF_i\|_{op} \leq e^{2(n-1)t_i}+C(n,\a) s_i^{1/2}\tau_i^{-1/2} = e^{2(n-1)t_i} + C(n,\alpha) s_i^{1/4} $;
\item[(B)] $|F_i(x)-\hat f_i(x)|\leq C(n)\tau_i^{1/2}+C(n,\a)s_i^{1/2} = C(n) s_i^{1/4} + C(n,\alpha)s_i^{1/2}$.
\end{enumerate}
We then define the truncated map $\tilde F_i:M^n\to\mathbb{R}^{N}$ given by 
\begin{equation}
\tilde F_i:\chi_i \hat f_i+(1-\chi_i) F_i.
\end{equation}

Since
\begin{equation}
d\tilde F_i= \chi_i d\hat f_i +(1-\chi_i) dF_i+ (\hat f_i-F_i) \cdot d\chi_i,
\end{equation}
we choose $s_i$ small enough so that $\|d\tilde F_i\|_{op}\leq 1+2\e_i$, and $F_i(M)\Subset U_i$ by (B), and $\hat f_i(x)\in \iota_i(\tilde M)$. We then claim that the map $f_i:M^n\to\tilde M$ given by $f_i:=\pi_i\circ \tilde F_i$ is the desired map. 

We first note that since $\pi_i\circ \hat f_i=\Phi$, the property $(\mathrm{III})$ follows directly from the construction, while $(\mathrm{I})$ follows from  $\|d\tilde F_i\|_{op}\leq 1+2\e_i$ and $\|d\pi_i\|_{op}\leq 1+\e_i$.

\medskip

It remains to show the conclusion $(\mathrm{II})$. 
It suffices to  note using Lemma~\ref{lma:d_0} and $\|d\pi_i\|_{op}\leq 1+\e_i$ that
\begin{equation}
\begin{split}
 \tilde d_{0}\left(f_i(x),\Phi(x)\right) &= \tilde d_{0}\left(\pi_i\circ \tilde F_i(x),\pi_i\circ \hat f_i(x)\right)\\&\leq d_{\tilde g(t_i)}\left(\pi_i\circ \tilde F_i(x),\pi_i\circ \hat f_i(x)\right)+C(n,\a)t_i^{1/2}\\
&\leq  (1+\e_i)|\tilde F_i(x)-\hat f_i(x)|+C(n,\a)t_i^{1/2}\\
&\leq C(n,\a) \left( s_i^{1/2}+s_i^{1/4}+t_i^{1/2}\right).
\end{split}
\end{equation}
This completes the proof.
\end{proof}

\subsection{Smoothing the initial isometry: the non-compact case}

In this subsection, we discuss how the initial distance isometry can be regularized in the non-compact case. For this purpose, we will use the standard Busemann function from  the work of Cheeger--Gromoll and Gromoll--Meyer \cites{CheegerGromoll1972,GromollMeyer1969}. 

We first recall its definition. Consider $(M^n,g)$ with $\mathrm{sec}\geq 0$, $x_0 \in M$ and the set $\mathcal{R}$ of geodesic rays $\gamma:[0,\infty)\rightarrow M$ with $\gamma(0)=x_0$. The standard Busemann function is defined as
\begin{equation}\label{eqn:standard-buse} 
b_{CG}(x):=\sup_{\gamma\in\mathcal{R}}\lim_{R\to+\infty} \left( R-d_g(x,\gamma(R))\right).
\end{equation}
It turns out to be more convenient for our purpose to consider the spherical Busemann function instead. Here, the spherical Busemann function $b:M\to \mathbb{R}$ with respect to $x_0$ is defined as
\begin{equation}\label{eqn:defn-buse}
    b(x):=\lim_{R\to+\infty} \left( R-d_g(x, \partial B_g(x_0,R))\right).
\end{equation}
By \cite[Theorem A]{Wu1979}, $b$ is $1$-Lipschitz and convex. Moreover, it is straightforward to check that $b$ dominates the standard Busemann and thus $b$ is also proper. 

Note that for the purpose of smoothing the initial isometry, i.e.~Proposition~\ref{prop:gauge-fixing-noncpt}, whether choosing $b$ to be the spherical Busemann function or the standard Busemann function $b_{CG}$ will not be relevant. It will only play a role in Lemma~\ref{lma:comp-flow-buse}.

\medskip

We will mimic the proof of Proposition~\ref{prop:gauge-fixing-cpt-1}. 
\begin{prop}\label{prop:gauge-fixing-noncpt}
Suppose the assumptions in Theorem~\ref{thm:compact-1} hold, and $M$ and $\tilde M$ are both non-compact. Then there exists $t_i\to0, 0<s_i<t_i$, $\Lambda_i\to 1$ and a sequence of smooth map $f_i:M\to \tilde M$ such that 
\begin{enumerate}
\item[(i)] $f_i^*\tilde g(t_i)\leq \Lambda_i g(s_i)$;
\item[(ii)] $f_i(M)\Subset \tilde M$, $\sup_M  |df_i|^2_{g(s_i),\tilde g(t_i)}+ |Ddf_i|_{g(s_i),\tilde g(t_i)}<+\infty$;
\item[(iii)] $\sup_M \tilde d_0( f_i(x),\Phi(x)) \to 0$ as $i\to+\infty$.
\end{enumerate}
\end{prop}
\begin{proof}
Let $R_i\to+\infty$, $t_i,\e_i\to 0$. Fix $x_0\in M$ and consider the Busemann function $b_i$ on $(\tilde M,g(t_i))$ with respect to $x_0\in M$. Since $b_i$ is proper in $(\tilde M,\tilde g(t_i))$, we might choose $L_i\to+\infty$ such that 
\begin{equation}
B_{\tilde g(t_i)}(x_0,R_i) \Subset \left\{ x\in M : b_i(x)\leq L_i\right\}.
\end{equation}
By \cite[Lemma 1]{GreeneShiohama}, there exists a continuous map $P_i:\tilde M \to \tilde M$ such that 
\begin{enumerate}
\item[(i)] $P_i$ is $1$-Lipschitz with respect to $\tilde g(t_i)$;
\item[(ii)] $P_i(\tilde M)=C_i:=\{ x\in \tilde M: b_i(x)\leq L_i\}$;
\item[(iii)] $P_i(z)=z$ if $z\in C_i$.
\end{enumerate}
We choose relatively compact smooth domains $V_i, W_i$ such that $C_i\subset W_i\Subset V_i\Subset \tilde M$. By suitably modifying the metric $\tilde g(t_i)$ near $\partial V_i$ and doubling $V_i$
, we might assume $C_i$ lives in a compact manifold. By Nash's embedding Theorem \cite[Theorem 3]{Nash}, there is a smooth embedding 
\begin{equation}
\iota_i: (W_i, \tilde g(t_i)) \to \mathbb{R}^{N}
\end{equation}
 for some $N=N(n)$.
We let $\pi_i$ be the nearest-point projection on an open neighbourhood $U_i$ of $\iota_i(W_i)\subseteq \mathbb{R}^{N}$ to $\iota_i(C_i)$ with $\|d\pi_i\|_{op}\leq 1+\e_i$. 

Now we consider the continuous map $\hat f_i: M\to \mathbb{R}^{N}$ given by 
\begin{equation}
\hat f_i:=\iota_i\circ P_i\circ \Phi
\end{equation}
which satisfies
\begin{equation}
\begin{split}
|\hat f_i(x)-\hat f_i(y)|&\leq d_{\tilde g(t_i)} \left(P_i\circ \Phi(x), P_i\circ \Phi(y) \right)\\
&\leq d_{\tilde g(t_i)} \left( \Phi(x),  \Phi(y) \right)\\
&\leq \tilde d_0(\Phi(x),\Phi(y)) =d_0(x,y)\\
&\leq d_{g(s)}(x,y)+C(n,\a)s^{1/2}
\end{split}
\end{equation}
for all $s\in (0,T]$. 

We choose $s_i$ small and  apply Proposition~\ref{prop:smoothing-map-Euc} to $(M,g(s_i))$ to obtain a smooth map $F_i:=F(\tau_i)$ for small $\tau_i=s_i^{1/2}$ so that $F_i: (M,g(s_i))\to \mathbb{R}^{N}$ is smooth and  satisfies 
\begin{enumerate}
\item[(A)] $\|dF_i\|_{op} \leq 1+C(n,\a) s_i^{1/2}\tau_i^{-1/2} = 1+C(n,\a) s_i^{1/4}$;
\item[(B)] $|F_i(x)-\hat f_i(x)|\leq C(n)s_i^{1/4}+C(n,\a)s_i^{1/2}$.
\end{enumerate}

Arguing as in the proof of Proposition~\ref{prop:gauge-fixing-cpt-1}, the map $f_i:=\pi_i\circ F_i$ is our desired map.
\end{proof}

\section{Uniqueness in the compact case with $\mathrm{sec}\geq 0$}

In this section, we will prove Theorem~\ref{thm:compact-1} when both manifolds are compact. We start by constructing a suitable harmonic map heat flow between the two flows. 
\begin{prop}\label{prop:Exist-HMHF-cpt-1}
Suppose the assumptions of Theorem~\ref{thm:compact-1} hold and suppose in addition that $M$ is compact. Then there exists a solution to the Ricci-harmonic map heat flow $F:(M,g(t))\to (\tilde M,\tilde g(t))$ for $t\in (0,T]$ such that $F$ is $1$-Lipschitz. Moreover $F(t)$ is homotopy equivalent to $\Phi$ and satisfies 
\begin{equation}
\sup_M \tilde d_0(F(x,t),\Phi(x))\leq C(n,\a)\sqrt{t}
\end{equation}
for all $t\in (0,T]$.
\end{prop}
\begin{proof}
By Proposition~\ref{prop:gauge-fixing-cpt-1} with $\Sigma=M$, there is $t_i,s_i\to 0$, $\Lambda_i\to 1$ and $f_i:M\to \tilde M$ such that 
\begin{enumerate}
\item[(I)] $f_i^*\tilde g(t_i)\leq \Lambda_i g(s_i)$ on $M$;
\item[(II)] $\tilde d_0(f_i(x),\Phi(x))\to 0$ as $i\to+\infty$. 
\end{enumerate} 
We consider the translated and rescaled Ricci flows
\begin{equation}
\left\{
\begin{array}{ll}
g^\sharp_i(t):= \Lambda_i g(s_i+\Lambda_i^{-1}t)\\[2mm]
\tilde g_i^\sharp(t):=\tilde g(t_i+t)
\end{array}
\right.
\end{equation}
for  $t\in [0,T)$ as $i\to+\infty$, such that $f_i$ is now $1$-Lipschitz with respect to $g_i^\sharp(0)$ and $\tilde g_i^\sharp(0)$.

We consider the Ricci-harmonic map heat flow: 
\begin{equation}
\left\{
\begin{array}{ll}
\left(\partial_t -\Delta_{g_i^\sharp(t), \tilde g_i^\sharp(t)} \right)F_i=0;\\[2mm]
F_i(0)=f_i.
\end{array}
\right.
\end{equation}
Since $g_i^\sharp$ and $\tilde g_i^\sharp$ both have bounded curvature (thanks to time translation) and $\mathrm{sec}\geq 0$, Lemma~\ref{lma:evo-eign} and the maximum principle imply that 
\begin{equation}
\left(F_i(t)\right)^* \tilde g_i^\sharp(t) \leq g_i^\sharp(t)
\end{equation}
on $M\times [0,T_{max})$. Parabolic regularity from \cite[Proposition 2.1]{Lee2025} allows us to let $i\to+\infty$ and obtain a limiting solution to the Ricci-harmonic map heat flow $F$ with respect to $g(t)$ and $\tilde g(t)$, on $M\times (0,T]$. Furthermore we have $|\partial_t F_i|\leq C(n,\a)t^{-1/2}$ uniformly for all $i\to+\infty$. By Lemma~\ref{lma:d_0}, 
\begin{equation}\label{eqn:F-close}
\begin{split}
\tilde d_0(F_i(x,t),\Phi(x)) &\leq \tilde d_0(F_i(x,t),f_i(x))+\tilde d_0(f_i(x,t),\Phi(x))\\
&\leq d_{\tilde g_i^\sharp(t)}(F_i(x,t),f_i(x))+C(n,\a)\sqrt{t_i+t}\\
&\quad+\tilde d_0(f_i(x,t),\Phi(x))\\
&\leq \int^{t}_0 |\partial_s F_i|_{\tilde g_i^\sharp (s)}\,ds +C(n,\alpha)\sqrt{t_i+t}\\
&\quad+\tilde d_0(f_i(x,t),\Phi(x))\\
&\leq C(n,\a)\sqrt{t_i+t} +\tilde d_0(f_i(x,t),\Phi(x)).
\end{split}
\end{equation}
Letting $i\to+\infty$ yields 
\begin{equation}
\sup_{ M}\tilde d_0(F(x,t),\Phi(x)) \leq C(n,\a) \sqrt{t}.
\end{equation}
In particular for each $t\in (0,T]$, $F(t)$ is homotopy equivalent to $\Phi$ and thus $\mathrm{deg}(F)=1$.
\end{proof}

Now we can prove the compact case of Theorem~\ref{thm:compact-1}. 
\begin{proof}[Proof of Theorem~\ref{thm:compact-1}, the compact case]

By Proposition~\ref{prop:Exist-HMHF-cpt-1}, there is a  smooth map $F:(M,g(t))\to (\tilde M,\tilde g(t))$  that is $1$-Lipschitz, of degree $1$ and satisfies 
\begin{equation}\label{eqn:dist-init-HMHF}
\tilde d_0(F(x,t),\Phi(x))\leq C(n,\a)\sqrt{t}.
\end{equation}
Therefore, 
\begin{equation}
\begin{split}
\mathrm{Vol}(\tilde M,\tilde g(t))&=\mathrm{deg}(F)\int_M d\mathrm{vol}_{F^*\tilde g(t)}\\
&=\int_M \sqrt{\frac{\det F^*\tilde g(t)}{\det g(t)}}\,d\mathrm{vol}_{g(t)}\\
&\leq \mathrm{Vol}(M,g(t)).
\end{split}
\end{equation}
Interchanging the role of $g(t)$ and $\tilde g(t)$ shows that the inequality is indeed equality. Hence $F(t)$ is a local isometry between $(M,g(t))$ and $(\tilde M,\tilde g(t))$. This also implies $\partial_t F=\Delta_{g(t),\tilde g(t)}F=0$ and hence $F$ is a constant map in time. By \eqref{eqn:dist-init-HMHF}, $F(x,t)\equiv \Phi(x)$ for all $x\in M$. This completes the proof.
\end{proof}

\section{Uniqueness in the compact case with regular sets}

In this section, we discuss how $\mathrm{sec}\geq 0$ can be relaxed under the assumption that the initial distance metric is Riemannian almost everywhere, i.e.~Theorem~\ref{thm:compact-2}. We first show that the Ricci flow is of almost $\mathrm{WPIC2}$ near the singular set.
\begin{prop}\label{prop:almostPIC2} Suppose the assumptions of Theorem~\ref{thm:compact-2} hold, and $x_0\in M$ and $r>0$ are such that $B_{d_0}(x_0,r)\Subset U$. Then for any $m\in\mathbb{N}$,  there is $\hat T_m>0$ independent of $x_0$, such that 
\begin{equation}
\mathrm{Rm}(g(x_0,t))+ t^m \mathrm{Id}\in \mathrm{C}_{\mathrm{PIC2}}
\end{equation}
for $t\in (0,\hat T_m]$.
\end{prop}
\begin{proof}

By \cite[Proposition 2.2]{BCW2019}, the function
\begin{equation}
v(x,t):=\inf\{ s\geq 0:  \Rm(x,t)+s\mathrm{Id} \in \mathrm{C}_{\mathrm{PIC2}}\} 
\end{equation}
satisfies 
\begin{equation}
\rheat v\leq \mathrm{scal}_{g(t)}\cdot v +C_1v^2
\end{equation}
in the sense of barrier, for some dimensional constant $C_1(n)>0$. 

By the maximum principle, for all $s<t<T$ the function $\hat v:=e^{-C_1t}v$ satisfies 
\begin{equation}
\begin{split}
\hat v(x_0,t)&\leq \int_M G(x_0,t;y,s) \hat v(y,s)\,d\mathrm{vol}_{g(y,s)},
\end{split}
\end{equation}
where $G(x,t;y,s)$ denotes the heat kernel of $\partial_t-\Delta_{g(t)}-\mathrm{scal}_{g(t)}$.

We split the integral into two parts: 
\begin{equation}
\begin{split}
&\quad \int_M G(x_0,t;y,s) \hat v(y,s)\,d\mathrm{vol}_{g(y,s)}\\
&=\left(\int_{B_{d_0}(x_0,r)} +\int_{M\setminus B_{d_0}(x_0,r)} \right) G(x_0,t;y,s) \hat v(y,s)\,d\mathrm{vol}_{g(y,s)}\\
&=:\mathbf{I}+\mathbf{II}.
\end{split}
\end{equation}

By Lemma~\ref{lma:d_0} and Colding's volume continuity under Gromov-Hausdorff convergence \cite{Colding}, we have $\mathrm{Vol}_{g(t)}\left(B_{g(t)}(x,\sqrt{t})\right)\geq v_0 t^{n/2}$ for all $x\in M$ and $t\in (0,T]$. Hence we might apply \cite[Proposition 3.1]{BCW2019} to see that  for fixed $t>0$ and as $s\to 0$ we have 
\begin{equation}
G(x,t;y,s)\leq \frac{C_2}{t^{n/2}}\cdot \exp\left(-\frac{d_{g(s)}^2(x,y)}{C_2t} \right)
\end{equation}
for some $C_2(n,\a,v_0)>0$.

Now we observe that for $y\notin B_{d_0}(x_0,r)$, we have 
\begin{equation}
d_{g(s)}(x_0,y)\geq d_0(x_0,y)-C(n,\a)\sqrt{s}\geq \frac12 r,
\end{equation}
for $s\to 0$, by Lemma~\ref{lma:d_0}. Putting this into $\mathbf{II}$ yields 
\begin{equation}
\begin{split}
\mathbf{II}&\leq \exp\left(-\frac{r^2}{8C_2t} \right)\cdot \int_{M\setminus B_{d_0}(x_0,r)} \frac{C_2}{t^{n/2}}\cdot \exp\left(-\frac{d_{g(s)}^2(x,y)}{2C_2t} \right) d\mathrm{vol}_{g(s)}\\
&\leq C_3 \exp\left(-\frac{r^2}{8C_3t} \right)
\end{split}
\end{equation}
as $s\to 0$. Here we have used that $\hat v\leq 1$ by assumption.

For $\mathbf{I}$, it follows from Lemma~\ref{lma:d_0} that $d\mathrm{vol}_{g(t)}\leq \mathcal{H}_{d_0}$ for all $t>0$. It follows from $\mathcal{H}^n_{d_0}(\Sigma)=0$, $\Rm(g_0)\in \mathrm{C}_{\mathrm{PIC2}}$, the dominated convergence theorem and heat kernel estimates that $\limsup_{s\to 0^+}\mathbf{I}\leq 0$. The result follows by choosing $\hat T_m$ small.
\end{proof}

\medskip
We prove an analogue of Proposition~\ref{prop:Exist-HMHF-cpt-1}.
\begin{prop}\label{prop:Exist-HMHF-cpt-2}
Under the assumptions in Theorem~\ref{thm:compact-2},  there exists a solution to the Ricci-harmonic map heat flow $F:(M,g(t))\to (\tilde M,\tilde g(t))$ for $t\in (0,T']$ for some $T'>0$ such that 
\begin{enumerate}
\item[(a)] $|dF|\leq 2n$;
\item[(b)] $
\sup_M \tilde d_0(F(x,t),x)\leq C(n,\a)\sqrt{t}$
\item[(c)] $F$ is smooth up to $t=0$ on $M_{reg}$;
\item[(d)] On a neighborhood $V$ of $\Sigma$, the following is true: for any $m\in \mathbb{N}$, there is $\hat T_m>0$ such that 
\begin{equation}
F^*\tilde g(t)\leq (1+t^m) g(t)
\end{equation}
on $V\times (0,\hat T_m]$.
\end{enumerate}
\end{prop}
\begin{proof}
By Proposition~\ref{prop:gauge-fixing-cpt-1}, there exist $t_i,s_i\to 0$, $\Lambda_i\searrow 1$ and a sequence of smooth maps $f_i:M\to M$ such that 
\begin{enumerate}
\item[(I)] $f_i^*\tilde g(t_i)\leq \Lambda_i g(s_i)$ on $M$;
\item[(II)] $\tilde d_0(f_i(x),x)\to 0$ as $i\to+\infty$;
\item[(III)] For any $\Omega\Subset M_{reg}$, there is $N$ such that $f_i|_{\Omega}=\mathrm{Id}$ for $i>N$. 
\end{enumerate}

We consider the translated rescaled Ricci flow
\begin{equation}
\left\{
\begin{array}{ll}
g^\sharp_i(t):= \Lambda_i g(s_i+\Lambda_i^{-1}t)\\[2mm]
\tilde g_i^\sharp(t):=\tilde g(t_i+t)
\end{array}
\right.
\end{equation}
for  $t\in [0,T)$ as $i\to+\infty$, such that $f_i$ is now $1$-Lipschitz with respect to $g_i^\sharp(0)$ and $\tilde g_i^\sharp(0)$. We consider the Ricci-harmonic map heat flow: 
\begin{equation}
\left\{
\begin{array}{ll}
\left(\partial_t-\Delta_{g_i^\sharp(t),\tilde g_i^\sharp(t)}\right)F_i=0;\\[2mm]
F_i(0)=f_i.
\end{array}
\right.
\end{equation}

We let $T_{max}>0$ be the maximal existence time of $F_i$ on $M\times [0,T_{max})$ such that $|dF_i|< 2n$. We claim that $T_{max}$ is bounded from below by a constant independent of $i\to+\infty$.

We denote $\Sigma_s:=\{x\in M: d_0(x,\Sigma)<s\}$ and let  $r_0>0$ be such that $\Sigma_{r_0}\Subset U$ and $B_{d_0}(x_0,r_0)\Subset U$ for all $x_0\in \Sigma_{4r_0}$. From Proposition~\ref{prop:almostPIC2}, we know that 
\begin{equation}\label{eqn:curv-almost}
\Rm(\tilde g(x_0,t))+t^m \mathrm{Id}\in \mathrm{C}_{\mathrm{PIC2}}
\end{equation}
for $t\in (0,\hat T_m]$, for any $x_0\in \Sigma_{4r_0}$.

Let $x\in \Sigma_{r_0}$. Then  Lemma~\ref{lma:d_0} implies
\begin{equation}
B_{g_i^\sharp(0)}(x,r_0)=B_{g(s_i)}(x, \Lambda_i^{-1/2} r_0)\Subset B_{d_0}(x,C_0\sqrt{s_i}+r_0)\subseteq B_{d_0}(x,2r_0)
\end{equation}
for $i$ sufficiently large. On the other hand, by \cite[Proposition 2.1]{Lee2025} and $|dF_i|<2n$ we have $|\partial_t F_i|\leq C(n,\a)t^{-1/2}$. If $y\in B_{g_i^\sharp(0)}(x,r_0)$, then same argument as in \eqref{eqn:F-close} implies 
\begin{equation}\label{eqn:F-close-i}
\begin{split}
d_{g_i^\sharp(0)}(F_i(y,t),y)&\leq d_{g_i^\sharp(0)}(F_i(y,t),f_i(y))+d_{g_i^\sharp(0)}(y,f_i(y))\\
&\leq C_0(n,\a)\sqrt{t}+o(1)
\end{split}
\end{equation}
so that $F_i(y,t)\in \Sigma_{4r_0}$ and hence \eqref{eqn:curv-almost} holds at $F_i(y,t)$, provided $i$ is sufficiently large and $t\leq \hat T(m,n,\a,r_0)\wedge T_{max}$.

By Lemma~\ref{lma:evo-eign}, \eqref{eqn:curv-almost} and $|dF_i|<2n$, the eigenvalue 
\begin{equation}
\varphi(x,t):=\sup\{ F_i^*\tilde g_i^\sharp(v,v): g_i^\sharp(v,v)=1\}
\end{equation}
satisfies 
\begin{equation}
\begin{split}
\left(\frac{\partial}{\partial t}-\Delta_{g_i^\sharp(t)} \right) \varphi&\leq 2n (\varphi-1)_+\cdot F_i^*\mathrm{scal}_{\tilde g_i^\sharp} +C_1(n)(t_i+t)^m 
\end{split}
\end{equation}
in the barrier sense, on $B_{g_i^\sharp(0)}(x,r_0)\times [0,\hat T_m\wedge T_{max}]$. Since $f_i$ is $1$-Lipschitz, $\varphi(0)\leq1$. Consider the function $v:=\varphi-1-(t_i+t)^\ell$. It satisfies $v(0)<0$ and 
\begin{equation}
\begin{split}
\left(\frac{\partial}{\partial t}-\Delta_{g_i^\sharp(t)} \right) v&\leq 2n \left(v+ (t_i+t)^\ell \right)\cdot F_i^*\mathrm{scal}_{\tilde g_i^\sharp} +C_1(t_i+t)^m-\ell (t_i+t)^{\ell-1}\\
&\leq 2n v\cdot  F_i^*\mathrm{scal}_{\tilde g_i^\sharp}
\end{split}
\end{equation}
in the barrier sense, on $B_{g_i^\sharp(0)}(x,r_0)\times [0,\hat T_m\wedge T_{max}]$, whenever $v>0$, provided that we choose $\ell_m:=2\cdot  \max\{m+1, C_1+2n\a\}$. We apply \cite[Theorem 1.1]{LeeTamCJ} (with scaling) to show that after shrinking $\hat T_m$ if necessary, 
\begin{equation}\label{eqn:varphi-esti-1}
F_i^*\tilde g_i^\sharp \leq \left( 1+(t_i+t)^{\ell_m} \right) g_i^\sharp
\end{equation}
on $\Sigma_{r_0}\times (0,\hat T_m\wedge T_{max}]$. 

We now estimate the Lipschitz constant outside $\Sigma_{r_0}$. Fix $x\notin \Sigma_{r_0}$ so that $B_{g_i^\sharp(0)}(x,r_0/3)\Subset M\setminus \Sigma_{r_0/2}$ by Lemma~\ref{lma:d_0}. Hence, \eqref{eqn:F-close-i} implies $F_i(y,t)\in M\setminus \Sigma_{r_0/4}$ for $y\in B_{g_i^\sharp(0)}(x,r_0/3)$ provided  $i$ is sufficiently large and $t\leq \hat T(m,n,\a,r_0)\wedge T_{max}$. By \cite[Corollary 3.2]{Chen2009}, $|\Rm(\tilde g_i^\sharp(t))|$ is uniformly bounded by some constant $\Lambda_0=\Lambda_0(n,\a,r_0,g_0)>1$ on $ M\setminus \Sigma_{r_0/4}$ for $t\in [0, T_{max})$ where $\tilde T=\tilde T(n,\a,r_0,g_0)>0$.

Hence the function $\varphi$ satisfies 
\begin{equation}
\begin{split}
\left(\frac{\partial}{\partial t}-\Delta_{g_i^\sharp(t)} \right) \varphi
&\leq \varphi (\varphi-1)_+ \cdot \left( F_i^*\mathrm{scal}_{\tilde g_i^\sharp}+n(n-1) \right)+2\varphi (n-1)\\
&\leq C_2(n)\Lambda_0 (\varphi-1)_+ +C_2(n)
\end{split}
\end{equation}
in the barrier sense, on $(M\setminus \Sigma_{r_0/4})\times [0,T_{max})$. Consider $u:=\varphi-1-\gamma C_2t$, which satisfies 
\begin{equation}
\begin{split}
\left(\frac{\partial}{\partial t}-\Delta_{g_i^\sharp(t)} \right)v &\leq C_2\Lambda_0 v+C_2+\gamma C_2^2\Lambda_0 t-\gamma  C_2 \\
&\leq C_2\left( 2n\Lambda_0 +1+\gamma C_2 \Lambda_0 t -\gamma \right) 
\end{split}
\end{equation}
whenever $v>0$. If we restrict $t\in [0, (2C_2\Lambda_0)^{-1}]$ and choose $\gamma:=2(1+2n\Lambda_0)$, then the maximum principle implies that its maximum on $M\setminus \Sigma_{r_0}$  must be attained on $\partial \Sigma_{r_0}$. Together with \eqref{eqn:varphi-esti-1}, we conclude that
\begin{equation}
\varphi\leq 1+\gamma C_2 t +(t_i+t)^{\ell_m}
\end{equation}
for $t\in [0,(2C_2\Lambda_0)^{-1}\wedge \hat T_m\wedge T_{max})$. In particular, $|dF|\leq 3n/2$ if $t$ is small enough. This proves a uniform lower bound for $T_{max}$ for all $i\to+\infty$. Hence \cite[Proposition 2.1]{Lee2025} enables us to take a convergent sub-sequence  $F_i\to F$ which is a smooth solution to the Ricci-harmonic map heat flow on $M\times (0,T']$ with $|dF|<2n$, for some $T'>0$. Furthermore, local parabolic theory implies that $F$ is smooth up to $t=0$ and equal to $F(0)=\mathrm{Id}$ on $M_{reg}\times [0,T']$.
\end{proof}

Now we are ready to complete the proof of Theorem~\ref{thm:compact-2}.
\begin{proof}[Proof of Theorem~\ref{thm:compact-2}]

It suffices to prove the uniqueness on $(0,S]$ for some $S>0$ by the uniqueness for smooth Ricci flows. We let $F$ be the harmonic map heat flow on $M\times (0,T']$, obtained from Proposition~\ref{prop:Exist-HMHF-cpt-2}. We also note that from Lemma~\ref{lma:d_0} and Colding's volume continuity \cite{Colding}, we might further assume 
\begin{enumerate}
\item $V_0^{-1}\leq \mathrm{Vol}(M,g(t))\leq V_0$;
\item $\inf_{x\in M}\mathrm{Vol}_{g(t)}\left(B_{g(t)}(x,1) \right)\geq v_0$,
\end{enumerate}
for some $V_0,v_0>0$ and all $t\in (0,T']$.

Since $F(0)=\mathrm{Id}$, $g(0)=\tilde g(0)=g_0$ on $M_{reg}$ smoothly, for any $m\in \mathbb{N}$ we might choose $\hat T_m$ such that outside the neighborhood of $\Sigma$, we have 
\begin{equation}\label{eqn:almost-1-lip}
F^*\tilde g(t)\leq (1+t^m) g(t)
\end{equation}
on $M\times (0,\hat T_m\wedge T']$. 

Since $F(t)$ is homotopy equivalent to the identity map, this implies 
\begin{equation}\label{eqn:vol-est-error}
\begin{split}
\mathrm{Vol}(\tilde M,\tilde g(t))&=\mathrm{deg}(F)\int_M d\mathrm{vol}_{\tilde g(t)}\\
&=\int_M \sqrt{\frac{\det F^*\tilde g(t)}{\det g(t)}}\,d\mathrm{vol}_{g(t)}\\
&\leq (1+t^m) ^{n/2}\cdot \mathrm{Vol}(M,g(t))\\
&\leq (1+nt^m)\cdot \mathrm{Vol}(M,g(t))\\
\end{split}
\end{equation}
for $t\in (0,T'\wedge \hat T_m]$, provided we shrink $\hat T_m$ further. Interchanging $g(t)$ and $\tilde g(t)$ yields 
\begin{equation}
(1+nt^m)^{-1}\cdot \mathrm{Vol}(M,g(t))\leq \mathrm{Vol}(\tilde M,\tilde g(t))\leq (1+nt^m)\cdot \mathrm{Vol}(M,g(t))
\end{equation}
for all $t\in (0,T'\wedge \hat T_m]$.

The function $u:=\frac1n \tr_g F^*\tilde g$ satisfies $u\leq 1+ t^m$ on $M\times (0,T'\wedge \hat T_m]$. Hence \eqref{eqn:vol-est-error} and the AM-GM inequality imply 
\begin{equation}
\begin{split}
\int_M \frac{u}{1+t^m}\,d\mathrm{vol}_{g(t)}&\geq 
\int_M \left(\frac{u}{1+t^m}\right)^{n/2}\,d\mathrm{vol}_{g(t)}\\
&\geq  (1+t^m)^{-n/2}\int_M \left(\frac{\det F^*\tilde g}{\det g}\right)^{1/2}\,d\mathrm{vol}_{g(t)}\\
&=(1+t^m)^{-n/2}\cdot \mathrm{Vol}(\tilde M,\tilde g(t))\\
&\geq (1+nt^m)^{-n/2-1}\cdot \mathrm{Vol}(M,g(t)).
\end{split}
\end{equation}
By replacing $m$ by $m+1$ and shrinking $\hat T_m$ if necessary, we conclude that 
\begin{equation}\label{eqn:esti-u}
\fint_M |u-1| d\mathrm{vol}_{g(t)}\leq t^m.
\end{equation}
By \cite[Proposition 2.1]{Lee2025}, $u$ satisfies 
\begin{equation}
|\nabla u|\leq C_n|dF||\nabla dF|\leq C_0(n,\a)t^{-1/2}.
\end{equation}
We fix an arbitrary $(x,t)\in M\times (0,T'\wedge \hat T_m]$ and consider $B_{g(t)}(x,\delta\sqrt{t})$ where 
\begin{equation}
u(z,t)\leq u(x,t)+C_0t^{-1/2}\cdot \delta t^{1/2}\leq u(x,t)+C_0\delta
\end{equation}
for $z\in B_{g(t)}(x,\delta\sqrt{t})$. 
Together with \eqref{eqn:esti-u}, this implies 
\begin{equation}
\begin{split}
1-u(x,t) &\leq C_0\delta+\fint_{B_{g(t)}(x,\delta \sqrt{t})} 1-u\,d\mathrm{vol}_{g(t)}\\
&\leq C_0\delta +\frac{\mathrm{Vol}_{g(t)}(M)}{\mathrm{Vol}_{g(t)}(x,\delta\sqrt{t})} t^m\\
&\leq C_0\delta+C_1(n,v_0,V_0) \delta^{-n} t^{m-n/2}.
\end{split}
\end{equation}
The argument is static. We are free to choose $\delta=t^\ell$ and $m=(n+1)(\ell+1)$. Since $\ell$ is arbitrarily large, after shrinking $\hat T_m$ further, we have shown that 
\begin{equation}
u(x,t)\geq 1-t^m
\end{equation}
on $M\times (0,T'\wedge \hat T_m]$. Combining this with \eqref{eqn:almost-1-lip} shows that
\begin{equation}
(1-t^m) g(t)\leq F^*\tilde g(t)\leq (1+t^m) g(t)
\end{equation}
on $M\times (0,T'\wedge \hat T_m]$. From this, we can argue as in \cite[Claim 3.3]{Lee2025} to show that $F:M\to \tilde M$ is a diffeomorphism, and hence $F_*g(t)$ is a Ricci-DeTurck flow with respect to $\tilde g(t)$ on $\tilde M\times (0,T'\wedge \hat T_m]$. Further, by volume comparison and \cite[Theorem 4.7]{CheegerGromovTaylor}, we have $\mathrm{inj}(\tilde{g}(t))\geq i_0\sqrt{t}$ for some $i_0>0$.

From now on, we freeze a large $m_0$.  By applying \cite[Theorem 4.1]{ChanLaiLee} to the time-translated Ricci flows, for some $\Lambda>0$ we have 
\begin{equation}
|F^*\tilde g(t)-g(t)|_{g(t)}\leq \Lambda |F^*\tilde g(s)-g(s)|_{g(s)} \leq \Lambda s^{m_0}
\end{equation}
for all $t\in [s,T'\wedge \hat T_{m_0}]$ and $s>0$. Here we remark that although \cite[Theorem 4.1]{ChanLaiLee} is stated for Ricci flow with $\mathrm{Rm}(\tilde g(t))\geq -1$, it only uses $\mathrm{sec}(\tilde g(t))\geq -1$. By letting $s\to0$, this shows that $F^*\tilde g(t)=g(t)$ and hence $\partial_t F\equiv 0$ on $M\times (0,T'\wedge \hat T_{m_0}]$. Since $F(0)=\mathrm{Id}$ on $M_{reg}$ and $\mathrm{Vol}_{g(t)}\left(\Sigma\right)\leq \mathcal{H}^n_{d_0}(\Sigma)=0$, $F(t)=\mathrm{Id}$ on $M\times (0,T'\wedge\hat T_{m_0}]$.  This proves uniqueness on $(0,T'\wedge \hat T_{m_0}]$, completing the proof.
\end{proof}

\section{Uniqueness in the non-compact case with $\mathrm{sec}\geq 0$}

In this section, we focus on the uniqueness problem in the non-compact case, i.e. the non-compact case in Theorem~\ref{thm:compact-1}. This relies on a weighted version of the proof in the compact case. We start by constructing a weight with good control. We will use the spherical Busemann function defined in \eqref{eqn:defn-buse}. 

Since uniqueness holds for Ricci flows with bounded curvature \cites{ChenZhu2006,Lee2025}, we will assume $T$ is small and might change from line to line.

\begin{lma}\label{lma:comp-flow-buse}
Suppose the assumptions in Theorem~\ref{thm:compact-1} hold and assume $M$ is non-compact. If $b_{g(t)}$ is the spherical Busemann function with respect to Ricci flow $g(t)$ for $t\in (0,T]$ and a fixed point $x_0$, then there exists $L_0>0$ such that for all $s\in (0,T]$ and $x\in M$,
\begin{equation}
b_{g(s)}(x)\geq \frac12 d_{g(s)}(x,x_0)-L_0.
\end{equation}
\end{lma}
\begin{proof}
Recall that from Lemma~\ref{lma:d_0}, we have
\begin{equation}
    d_{g(t)}\leq d_{g(s)}\leq d_{g(t)}+C_0(n,\a)\sqrt{t-s}
\end{equation}
for all $0<s<t<T$. This implies 
\begin{equation}
B_{g(s)}(x_0,R)\subseteq B_{g(t)}(x_0,R)\subseteq B_{g(s)}(x_0,R+C_0\sqrt{t-s})
\end{equation}
for all $s<t$ and $R\to+\infty$. In particular, 
\begin{equation}
\begin{split}
d_{g(s)}\left(x,\partial B_{g(s)}(x_0,R)\right)&\leq 
d_{g(s)}\left(x,\partial B_{g(t)}(x_0,R)\right)\\
&\leq d_{g(t)}\left(x,\partial B_{g(t)}(x_0,R)\right)+C_0\sqrt{t-s}
\end{split}
\end{equation}
and hence, 
\begin{equation}\label{eqn:sph-buse}
\begin{split}
b_{g(s)}(x)&\geq b_{g(t)}(x)-C_0\sqrt{t-s}
\end{split}
\end{equation}
for $s<t$. By taking $t=T$, we have shown that 
\begin{equation}
b_{g(s)}(x)\geq b_{g(T)}(x)-C_1(n,\a,T)
\end{equation}

Furthermore, it follows from the definition that $b_{g(T)}$ dominates the Busemann function $b_{CG}$ considered by Cheeger-Gromoll \cite{CheegerGromoll1972} (see \eqref{eqn:standard-buse}) which is known to be asymptotic to $d_{g(T)}(x,x_0)$ when $x\to +\infty$, see \cite[Lemma 1 (b)]{Drees}. Hence, 
\begin{equation}
b_{g(s)}(x)\geq \frac12 d_{g(T)}(x,x_0)-C_2
\end{equation}
for some large constant $C_2>0$. The result follows from applying Lemma~\ref{lma:d_0} to compare $d_{g(T)}$ and $d_{g(s)}$.
\end{proof}

We aim to evolve the spherical Busemann function along the Ricci flow. In the following, for a smooth function $u:M\times I\to \R$, where $I\subset \R$ is an interval, and $(g(t))_{t\in I}$ a Ricci flow on $M$, we say that $u$ solves the heat equation coupled with Ricci flow, provided
\begin{equation}
    \left(\frac{\partial}{\partial t} - \Delta_{g(t)}\right)u=0\, .
\end{equation}

The next Lemma will be used to show that the heat flow preserves the convexity. 
\begin{lma}\label{lma:evo-convx}
Assume $u$ solves the heat equation coupled with Ricci flow. Assume $\mathrm{sec}(g(t))\geq 0$ for all $t\in I$. Then the function given by
\begin{equation}
\varphi(x,t):=\inf\{s\geq 0: \nabla^2 u+ s g\geq 0\}
\end{equation}
satisfies 
\begin{equation}
\rheat \varphi\leq \mathrm{scal}(g(t))\cdot \varphi
\end{equation}
in the barrier sense, whenever $\varphi>0$.
\end{lma}
\begin{proof}
The proof follows that of \cite[Lemma 3.2]{ChanLeePeachey}. For completeness we outline the proof below.
From \cite[Corollary  2.34]{ChowBOOK}, $u_{ij}$ satisfies 
\begin{equation}
\left(    \partial_t-\Delta_{g(t)}\right)u_{ij}=2R_{iklj} u_{kl}-R_{ik} u_{kj}-R_{jk}u_{ik},
\end{equation}
in orthonormal frame of $g(t)$. In particular, after applying the vector bundle trick as in the proof of Lemma~\ref{lma:evo-eign}, at $(x_0,t_0)$ where $\varphi(x_0,t_0)>0$, we can assume $\varphi(x_0,t_0)=-u(E_1,E_1)(x_0,t_0)$ for some orthonormal frame $\{E_i\}_{i=1}^n$ defined on a space-time neighbourhood of $(x_0,t_0)$, where $E_1$ is the eigenvector of $\nabla^2 u$ at $(x_0,t_0)$, and
\begin{equation}
    \begin{split}
        \left(    \partial_t-\Delta_{g(t)}\right)  (-u(E_1,E_1))&\leq -2R_{1kk1} u_{kk}\\
        &\leq 2R_{1kk1}\varphi \\
        &=2R_{11}\varphi\\
        &\leq \mathrm{scal}\cdot \varphi
    \end{split}
\end{equation}
where we have used \cite[Lemma 3.2]{ChanLee2025}.
\end{proof}

We now construct a proper function in the Ricci flow background which is convex and of linear growth.
\begin{prop}\label{prop:convex-prop}
Suppose the assumptions in Theorem~\ref{thm:compact-1} hold and assume $M$ is non-compact. Then for $x_0\in M$ there is a smooth solution $u_\infty$ to the heat equation coupled with Ricci flow on $M\times (0,T]$ and constants $L_1,L_2 \in [0,\infty)$ such that 
\begin{enumerate}
\item for all $(x,t)\in M\times (0,T]$,
\begin{equation}
d_0(x,x_0)+L_2\sqrt{t}\geq u_\infty(x,t)\geq \frac12 d_{0}(x,x_0)-L_1-L_2\sqrt{t};
\end{equation}
\item $|\nabla u_\infty|\leq 1$ on $M\times (0,T]$;
\item $\nabla^2 u_\infty\geq 0$ on $M\times (0,T]$.
\end{enumerate}
\end{prop}
\begin{proof}
Take sequences $\e_i,s_i\in(0,1)$ with $\e_i,s_i\to 0$, and consider the metric $g(t_i)$ on $M$. We let $b_{g(t_i)}$ be the spherical Busemann function with respect to $g(t_i)$ and a fixed point $x_0\in M$. By \cite[Proposition 2.3]{GreeneWu}, we might find a sequence of smooth functions $\b_i$ on $M$ such that on $M$, 
\begin{enumerate}
\item[(i)] $\nabla^{2,g(t_i)}\b_i\geq -\e_ig(t_i)$; 
\item[(ii)] $|\nabla^{g(t_i)}\b_i|\leq 1+\e_i$;
\item[(iii)] $|\b_i-b_{g(t_i)}|\leq \e_i$.
\end{enumerate}

We let $K(x,t;y,s)$ be the heat kernel with respect to the Ricci flow $g(t)$ on $M\times (0,T]$. We let 
\begin{equation}
u_i(x,t):=\int_M K(x,t;y,t_i)\cdot \b_i(y)\,d\mathrm{vol}_{g(t_i)}
\end{equation}
on $M\times [t_i,T]$. That is the minimal solution to the heat equation coupled with $g(t)$, starting from $u(t_i)=\b_i$.

\begin{claim}\label{claim:zero}
There exists $C_1(n,\a)>0$ such that for all $i\to+\infty$, $(x,t)\in M\times (t_i,T]$, 
\begin{equation}
|u_i(x,t)|\leq C_1\left(d_{g(t_i)}(x,x_0)+\e_i+\sqrt{t} \right).
\end{equation}
In particular, $u_i$ sub-converges to a solution $u_\infty$ to the heat equation on $M\times (0,T]$ such that 
\begin{equation}
|u_\infty(x,t)|\leq C_1\left(d_{0}(x,x_0)+\sqrt{t} \right).
\end{equation}
\end{claim}
\begin{proof}
By \cite[Lemma 3.3]{ChanLee2025} (see also \cite[Proposition 3.1]{BCW2019}), $\mathrm{scal}(g(t))\geq 0$ and the maximum principle,
\begin{equation}\label{eqn:Heat-ker}
K(x,t;y,t_i)\leq \frac{C_0}{\mathrm{Vol}_{g(t_i)}(x,\sqrt{t})}\cdot \exp\left( -\frac{d_{g(t_i)}(x,y)^2}{C_0t}\right)
\end{equation}
for all $x,y\in M$ and $2t_i< t\leq T$. Here $C_0=C_0(n,\a)>0$. Hence for any $x\in M$,  we have 
\begin{equation}
\begin{split}
|\b_i(y)|&\leq |\b_i(y)-\b_i(x)|+|\b_i(x)-\b_i(x_0)|+|\b_i(x_0)-b_{g(t_i)}(x_0)|\\
&\leq (1+\e_i) d_{g(t_i)}(x,y)+(1+\e_i)d_{g(t_i)}(x,x_0)+\e_i 
\end{split}
\end{equation}
 and thus
\begin{equation}
\begin{split}
|u_i&(x,t)|\\
&\leq \int_M \frac{C_0}{\mathrm{Vol}_{g(t_i)}(x,\sqrt{t})} \exp\left( -\frac{d_{g(t_i)}(x,y)^2}{C_0t}\right) |\b_i(y)|\,d\mathrm{vol}_{g(t_i)}\\
&\leq 2\int_M \frac{C_0}{\mathrm{Vol}_{g(t_i)}(x,\sqrt{t})} \exp\left( -\frac{d_{g(t_i)}(x,y)^2}{C_0t}\right) d_{g(t_i)}(x,y)\,d\mathrm{vol}_{g(t_i)}\\
&\quad +\left(2d_{g(t_i)}(x,x_0)+\e_i \right)\int_M \frac{C_0}{\mathrm{Vol}_{g(t_i)}(x,\sqrt{t})} \exp\left( -\frac{d_{g(t_i)}(x,y)^2}{C_0t}\right) \,d\mathrm{vol}_{g(t_i)}\\
&=:\mathbf{I}+\mathbf{II}.
\end{split}
\end{equation}

We now estimate the the integrals $\mathbf{I}$ and $\mathbf{II}$. We only estimate $\mathbf{I}$, since $\mathbf{II}$ is similar and simpler. By the co-area formula, integration by parts and volume comparison, 
\begin{equation}
    \begin{split}
        &\quad \int_M \frac{C_0}{\mathrm{Vol}_{g(t_i)}(x,\sqrt{t})}\cdot \exp\left( -\frac{d_{g(t_i)}(x,y)^2}{C_0t}\right) d_{g(t_i)}(x,y)\,d\mathrm{vol}_{g(t_i)}\\
        &=\int^{+\infty}_0 \frac{C_0r}{\mathrm{Vol}_{g(t_i)}(x,\sqrt{t})}\cdot \exp\left( -\frac{r^2}{C_0t}\right) \cdot |\partial B_{g(t_i)}(x,r)| \,dr\\
        &\leq \int^{+\infty}_0 \frac{\mathrm{Vol}_{g(t_i)}(x,r)}{\mathrm{Vol}_{g(t_i)}(x,\sqrt{t})} \cdot \frac{2r^2}{t} \exp\left(-\frac{r^2}{C_0t}\right)dr\\
        &=\left(\int^{\sqrt{t}}_0+\int^{+\infty}_{\sqrt{t}} \right) \frac{\mathrm{Vol}_{g(t_i)}(x,r)}{\mathrm{Vol}_{g(t_i)}(x,\sqrt{t})} \cdot \frac{2r^2}{t}  \exp\left(-\frac{r^2}{C_0t}\right)dr\\
        &\leq \int^{\sqrt{t}}_0  \frac{2r^2}{t} \exp\left(-\frac{r^2}{C_0t}\right)dr +\int^{+\infty}_{\sqrt{t}}  \left( \frac{r}{\sqrt{t}}\right)^n \cdot \frac{2r^2}{t}  \exp\left(-\frac{r^2}{C_0t}\right)dr\\
        &\leq C_2 \sqrt{t}.
    \end{split}
\end{equation}
To summarize, we have shown that 
\begin{equation}
    |u_i(x,t)|\leq C_3\left(d_{g(t_i)}(x,x_0)+\e_i+\sqrt{t} \right).
\end{equation}
The sub-sequential convergence in the $C^\infty_{loc}$ topology follows from standard parabolic theory.
\end{proof}

We next show that $u_\infty$ is $1$-Lipschitz and convex for each $t\in (0,T]$.
\begin{claim}\label{claim:lip-conve}
We have $|\nabla u_\infty|\leq 1$ and $\nabla^2 u_\infty\geq 0$ on $M\times (0,T]$.
\end{claim}
\begin{proof}

We make use of the construction of $u_\infty$. It is standard that along the Ricci flow we have
\begin{equation}\label{eqn:evo-Lip-flow}
\rheat |\nabla u_i|^2=-2|\nabla^2 u_i|^2\leq 0,
\end{equation}
for example see \cite[Lemma 2.3]{Lee2025}. The main difficulty is to justify the maximum principle. To do this, we fix any $r>0$ and consider $B_{g(t_i)}(x_0,r)\times [t_i,T]$. For large $r>1$, we consider the rescaled solution $\hat u_i(x,t):=r^{-1}u_i(t_i+r^2 t)$ with respect to $\hat g_i(t):=r^{-2}g(t_i+r^2t)$. It satisfies $|\hat u_i|\leq C_2(n,\a,T)$ on $B_{\hat g_i(0)}(x_0,2)\times [0,r^{-2}(T-t_i)]$ by  Claim~\ref{claim:zero}. By an argument similar to the proof of \cite[Proposition 2.1]{Lee2025}, we might assume 
\begin{equation}\label{eqn:ap-C12}
|\hat\nabla \hat u_i|^2+|\hat \nabla^2 \hat u_i|\leq C_3(n,\a,T)t^{-1}
\end{equation}
on $B_{\hat g_i(0)}(x_0,1)\times (0,r^{-2}(T-t_i)]$ if $r$ is large enough.

Since the initial data satisfies $|\nabla^{g(t_i)}\b_i|\leq 1+\e_i$, \eqref{eqn:evo-Lip-flow} and \eqref{eqn:ap-C12} enable us to apply \cite[Theorem 1.1]{LeeTamCJ} and show that 
\begin{equation}\label{eqn:grad-appr}
|\nabla \hat u_i(x,t)|\leq 1+\e_i+t
\end{equation}
on $B_{\hat g_i(0)}(x_0,1/2)\times (0,r^{-2}(T-t_i)\wedge \hat T]$, for some $\hat T(n,\a,T)>0$. The gradient estimate follows by scaling it back and letting $r\to+\infty$, using Lemma~\ref{lma:d_0}.

It remains to show convexity. Since the initial data satisfies $\nabla^{2,g(t_i)}\b_i\geq -\e_i g(t_i)$, $\hat u_i$ satisfies 
\begin{equation}
\hat\nabla^2 \hat u_i(0)\geq -\e_i r \hat g_i(0).
\end{equation}

We define the barrier function 
\begin{equation}
\hat v_i(x,t):=\int_M \hat G_i(x,t;y,0) \cdot (\e_i r)\,d\mathrm{vol}_{\hat g_i(0)}
\end{equation}
so that the eigenvalue 
\begin{equation}
\hat \varphi_i(x,t):=\inf\{ s\geq 0 : \hat\nabla^2 \hat u_i+s \hat g_i \geq 0\}
\end{equation}
satisfies $\hat\varphi_i- \hat v_i \leq 0$ at $t=0$. Now \eqref{eqn:ap-C12} and Lemma~\ref{lma:evo-convx} allow us to apply \cite[Theorem 1.1]{LeeTamCJ} to conclude that 
\begin{equation}
\hat \varphi_i(x,t)\leq \hat v_i(x,t)+t\leq C_4 \e_i r+t
\end{equation}
on $B_{\hat g_i(0)}(x_0,1/2)\times (0,r^{-2}(T-t_i)\wedge \hat T]$, where we have used the heat kernel estimate from  \cite[Lemma 3.3]{ChanLee2025}. By scaling it back, letting $r\to+\infty$ and applying Lemma~\ref{lma:d_0}, we have shown that $\nabla^2 u_\infty\geq 0$ on $M\times (0,T]$.
\end{proof}

The gradient estimate in particular allow us to improve the regularity so that $u_\infty$ is comparable to distance function. 
\begin{claim}\label{claim:proper}
The function $u_\infty$ is proper and satisfies 
\begin{equation}
d_0(x,x_0)+L_2\sqrt{t}\geq u_\infty(x,t)\geq \frac12 d_{0}(x,x_0)-L_1-L_2\sqrt{t}
\end{equation}
on $M\times (0,T]$, for some $L_1,L_2>0$.
\end{claim}
\begin{proof}
We use the estimate along the approximation. By \eqref{eqn:grad-appr} and \cite[Proposition 2.1]{Lee2025},  we have 
\begin{equation}
|\hat \nabla^2 \hat u_i(x,s)|\leq C_5 s^{-1/2}
\end{equation}
on $B_{\hat g_i(0)}(x_0,1/2)\times (0,r^{-2}(T-t_i)\wedge \hat T]$ for some $\hat T>0$. We integrate it from $s=0$ to $s=t$ to conclude that 
\begin{equation}
|\hat u_i(x,t)-\hat u_i(x,0)|\leq C_6 \sqrt{t}.
\end{equation}
Scaling back, this translates to
\begin{equation}\label{eqn:u_i-initial}
|u_i(x,t)-\b_i(x)|\leq C_6\sqrt{t-t_i}
\end{equation}
for all $(x,t)\in M\times [t_i,T]$. Combining this with Lemma~\ref{lma:comp-flow-buse} gives 
\begin{equation}
\begin{split}
u_i(x,t)&\geq \b_i(x)-C_6\sqrt{t-t_i}\\
&\geq b_{g(t_i)}(x)-\e_i-C_6\sqrt{t-t_i} \\
&\geq \frac12 d_{g(t_i)}(x,x_0)-L_0-\e_i-C_6\sqrt{t-t_i},
\end{split}
\end{equation}
while 
\begin{equation}
\begin{split}
u_i(x,t)&\leq \b_i(x)+C_6\sqrt{t-t_i}\\
&\leq b_{g(t_i)}(x)-b_{g(t_i)}(x_0)+C_6\sqrt{t-t_i}+\e_i\\
&\leq (1+\e_i)d_{g(t_i)}(x,x_0)+C_6\sqrt{t-t_i}+\e_i.
\end{split}
\end{equation}
This completes the proof of Claim \ref{claim:proper} by letting $i\to+\infty$ and using Lemma~\ref{lma:d_0}.
\end{proof}
The proof is completed by combining all claims.
\end{proof}

\medskip

We now construct a harmonic map heat flow between two Ricci flows. 
\begin{prop}\label{prop:Exist-HMHF-noncpt}
Under the assumptions of Theorem~\ref{thm:compact-1}, supposing in addition that $M$ is non-compact, there exists a solution to the Ricci-harmonic map heat flow $F:(M,g(t))\to (\tilde M,\tilde g(t))$ for $t\in (0,T]$ such that $F$ is $1$-Lipschitz. Moreover $F(t)$ is homotopy equivalent to $\Phi$ and satisfies 
\begin{equation}\label{eqn:initial-F-noncpt}
\sup_M \tilde d_0(F(x,t),\Phi(x))\leq C(n,\a)\sqrt{t}
\end{equation}
for all $t\in (0,T]$.
\end{prop}
\begin{proof}
The proof is identical to that of Proposition~\ref{prop:Exist-HMHF-cpt-1}, except we use Proposition~\ref{prop:gauge-fixing-noncpt} and the short-time existence theory of the harmonic map heat flow from \cites{ChenZhu2006,Lee2025}.
\end{proof}

Now we are ready to prove the uniqueness theorem, i.e. Theorem~\ref{thm:compact-1}, in the non-compact case. 

\begin{proof}[Proof of Theorem~\ref{thm:compact-1}, the non-compact case]

Let $x_0\in M$ and $\tilde x_0:=\Phi(x_0)\in \tilde M$. We let $u_\infty$ and $\tilde u_\infty$ be the convex proper functions obtained from applying Proposition~\ref{prop:convex-prop} to $\left(M,g(t),x_0\right)$ and $(\tilde M,\tilde g(t),\tilde x_0)$ respectively. We also let $F:M\to\tilde M$ be the Ricci-harmonic map heat flow obtained from Proposition~\ref{prop:Exist-HMHF-noncpt}.

\begin{claim}\label{claim:weigh}
The function $f:=F^*\tilde u_\infty-u_\infty$ is non-negative on $M\times (0,T]$.
\end{claim}
\begin{proof}
For convenience, we use the construction of both $u_\infty$ and $\tilde u_\infty$. The function $f$ satisfies 
\begin{equation}
\begin{split}
\rheat f&=\rheat F^*\tilde u_\infty\\
&=\left(\tilde g^{\a\b} -g^{ij} F^\a_i F^\b_j\right) \tilde\nabla_\a \tilde\nabla_\b \tilde u_\infty \geq 0
\end{split}
\end{equation}
on $M\times (0,T]$, where we have used the convexity of $\tilde u_\infty$ and the fact that $F$ is $1$-Lipschitz.

For convenience, we trace back the construction of $u_\infty$ and $\tilde u_\infty$. From \eqref{eqn:u_i-initial}, we have both
\begin{equation}
\left\{
\begin{array}{ll}
\displaystyle |u_i(x,t)-b_{g(t_i)}(x)|\leq C_0\sqrt{t-t_i}\\[2mm]
\displaystyle|\tilde u_i(x,t)-b_{\tilde g(t_i)}(x)|\leq C_0\sqrt{t-t_i}
\end{array}
\right.
\end{equation}
while from Lemma~\ref{lma:d_0}, $\tilde x_0=\Phi(x_0)$ and  the derivation of \eqref{eqn:sph-buse}, we also have 
\begin{equation}
|b_{g(t_i)}(x)-\Phi^* b_{\tilde g(t_i)}(x)|\leq C_1\sqrt{t-t_i}.
\end{equation}

This together with \eqref{eqn:initial-F-noncpt}, Lemma~\ref{lma:d_0} and $|\tilde \nabla \tilde u_\infty|\leq 1$ shows that 
\begin{equation}
\begin{split}
&\quad |u_\infty(x,t)-\tilde u_\infty\left( F(x,t),t\right)|\\
&\leq |u_\infty(x,t)-\tilde u_\infty\left( \Phi(x),t\right)|+|\tilde u_\infty(\Phi(x),t)-\tilde u_\infty\left( F(x,t),t\right)|\\
&\leq |u_\infty(x,t)-\tilde u_\infty\left( \Phi(x),t\right)|+d_{\tilde g(t)}\left(\Phi(x),F(x,t) \right)\\
&\leq C_2\sqrt{t}
\end{split}
\end{equation}
for $(x,t)\in M\times (0,T]$.

 This shows that the function $f$ is bounded on $M\times (0,T]$, and satisfies $f(t)=O(t^{1/2})$ as $t\to 0$. Since $g(t)$ has bounded curvature for $t>0$, by applying the minimum principle to $f$ on $M\times [s,T]$ and letting $s\to 0$, we conclude that $f(x,t)\geq 0$ on $M\times (0,T]$.
\end{proof}

Thanks to the linear growth from Proposition~\ref{prop:convex-prop} and volume comparison, 
\begin{equation}
E(t):=\int_M e^{-u_\infty}\,d\mathrm{vol}_{g(t)}
\end{equation}
is a well-defined energy for $t\in (0,T]$. Moreover from Claim~\ref{claim:weigh}, the fact that $F$ is $1$-Lipschitz, homotopy equivalent to the isometry $\Phi$ and proper by \eqref{eqn:initial-F-noncpt},
\begin{equation}
\begin{split}
E(t)&=\int_M e^{-u_\infty}\,d\mathrm{vol}_{g(t)}\\
&\geq \int_M e^{-F^*\tilde u_\infty}\sqrt{\frac{\det F^*\tilde g}{\det g}}\,d\mathrm{vol}_{g(t)}\\
&=\int_{\tilde M} e^{-\tilde u_\infty}\,d\mathrm{vol}_{\tilde g(t)}=:\tilde E(t).
\end{split}
\end{equation}
By interchanging the role of $\big(M,g(t),u_\infty(t)\big)$ and $\big( \tilde M,\tilde g(t),\tilde u_\infty(t)\big)$, we conclude that $E(t)\equiv \tilde E(t)$ for all $t\in (0,T]$. This forces the inequalities to be equalities. In particular, $F^*\tilde g(t)\equiv g(t)$ on $M\times (0,T]$. This shows that $F$ is static in time and completes the proof by \eqref{eqn:initial-F-noncpt}. 
\end{proof}

\end{document}